\documentclass[11pt]{amsart}
\usepackage{amssymb,amsmath,txfonts,amsrefs, amscd}
\usepackage{appendix}
\usepackage[all]{xy}
\newtheorem{thm}{Theorem}[section]
\newtheorem{prop}[thm]{Proposition}
\newtheorem{lemma}[thm]{Lemma}
\newtheorem{remark}[thm]{Remark}

\newtheorem{definition}[thm]{Definition}
\newtheorem{cor}[thm]{Corollary}

\numberwithin{equation}{section}

\def\pa{\partial}
\def\bpa{\bar{\partial}}
\def\dpp{{\partial_+}}
\def\dpm{{\partial_-}}
\def\dpps{{\partial^*_+ }}
\def\dpms{{\partial^*_-}}  
\def\om{\omega}
\def\w{\wedge}
\def\dpmm{{\partial'_-}}
\def\Dmm{{\Delta'_-}}

\def\cok{{\rm coker}}

\begin{document}

\title[]{Hodge Theory and Symplectic Boundary Conditions}
\author{ Li-Sheng Tseng}
\address{Department of Mathematics, University of California, Irvine 92697 CA}
\email{lstseng@math.uci.edu}
\author{Lihan Wang}
\address{Department of Mathematics, University of California, Irvine 92697 CA}
\email{lihanw@uci.edu}

\date{\today}


\begin{abstract}
We study symplectic Laplacians on compact symplectic manifolds with boundary. These Laplacians are associated with symplectic cohomologies of differential forms and can be of fourth-order.  We introduce several natural boundary conditions on differential forms and use them to establish Hodge theory by proving various form decomposition and also isomorphisms between the symplectic cohomologies and the spaces of harmonic fields.  These novel boundary conditions can be applied in certain cases to study relative symplectic cohomologies and Lefschetz maps between relative de Rham cohomologies. As an application, our results are used to solve boundary value problems of differential forms. 


\end{abstract}

\maketitle
\tableofcontents

\section{Introduction}

On a symplectic manifold $(M^{2n}, \om)$, there is a natural decomposition of the standard exterior derivative operator \cite{TY2}
$$ d = \dpp + \om \w \dpm~.$$  
The pair $(\dpp, \dpm)$ are dependent on the symplectic structure $\om$ and are linear first-order differential operators with properties similar to the familiar Dolbeault operators $(\partial, \bar{\partial})$ of complex geometry.  Of importance, these operators are suggestive of a new type of analysis on symplectic manifolds.  For $(\dpp, \dpm)$ should be thought of as the fundamental building blocks to write down linear elliptic operators that are inherently symplectic.  And the analysis of the global spectral properties of such elliptic operators should result in interesting symplectic invariants.

In this paper, we shall mainly study four symplectic elliptic operators.   Recall that on any symplectic manifold, there exists a compatible triple, $(\omega, J, g)$, involving additionally an almost complex structure $J$, and a Riemannian metric $g$.  Using the standard definition of the Riemannian inner product and adjoint for operators, we define the following elliptic operators
\begin{align*}
\Delta_{+}&= \dpps\dpp+\dpp\dpps,~~\text{ on $P^{k}$ with $k<n$},\\ 
\Delta_{-}&=\dpms\dpm+\dpm\dpms,~~ \text{ on $P^{k}$ with $k<n$},\\
\Delta_{++}&=(\dpp\dpm)^{\ast}\dpp\dpm+(\dpp\dpps)^2,~~ \text{ on $P^{n}$ },\\
\Delta_{--}&=\dpp\dpm(\dpp\dpm)^{\ast}+(\dpms\dpm)^2,~~ \text{ on $P^{n}$  }
\end{align*}
These second- and fourth-order operators act on the space of primitive forms $P^k$, where $k=1,\ldots, n$.  Primitive forms can be heuristically thought of as forms that are orthgonal to $\om$ and are analogus to the holomorphic forms of complex geometry.   
We will call the above operators symplectic Laplacians as  they are the Laplacian operators associated with the symplectic cohomologies of differential forms discussed in \cite{TY2}.  We do point out that our definition of the fourth-order Laplacians are different from those in \cite{TY1, TY2} as our modifications ensure that $\Delta_{++}$ and $\Delta_{--}$ are elliptic on $P^n$.

For these symplectic Laplacians, we begin the study of their spectral properties in this paper by analyzing their Hodge theory on symplectic manifolds with boundary.  We here emphasize Hodge theory as it concerns the zero modes of the Laplacians and is also the basic tool with wide-ranging applications in the study of manifolds with boundary.   A concurrent motivation for us is to use Hodge theory to study the symplectic cohomologies of forms in \cite{TY1, TY2} on compact manifolds with boundary.

Unlike the case of closed manifolds, statements of Hodge theory in the case of compact manifolds with boundary are more subtle and requires more than just the ellipticity of the Laplacian operator of interest.  Boundary conditions must be placed on differential forms and sometimes also on the boundary of the manifold in order for Hodge theory to work.  For example, we recall that Friedrichs \cite{F} and Morrey \cite{M1} extended the classical Hodge theory for the Laplace-Beltrami operator $\Delta_d=d\, d^{\ast} + d^{\ast} d$ to the case of manifold with boundary.  They studied the subspaces of harmonic fields which are forms that are both $d$-closed and $d^{\ast}$-closed.  (Note the distinction in the boundary case: a harmonic {\it form}, that is a zero of the Laplacian, is not necessarily a harmonic {\it field}.) They showed that the space of harmonic fields satisfying either the Dirichlet (denoted here by $D$) or Neumann (denoted by $N$) boundary condition is finite-dimensional and that several types of decompositions of differential forms hold with respect to these two boundary conditions.  In contrast, for the $\bar{\partial}-$Laplace operator $\Delta_{\bar{\partial}}=\bpa\, \bpa^{\ast}+ \bpa^{\ast}\bpa$ on complex manifolds with boundary, a consistent Hodge theory requires the boundary to satisfy the strong pseudo-convex condition. Moreover, the boundary condition on differential forms is the $\bar{\partial}-$Neumann condition. In Table 1, we summarize the well-known boundary conditions involved in Hodge theory for these two cases, with $\rho$ denoting a boundary defining function. (For a general reference, see \cite{G} for the Riemannian case and \cite{M2} for the complex case.) 

\begin{table}[!t]
\caption{Standard boundary conditions for Hodge theory on manifolds with boundary. }
\renewcommand{\arraystretch}{1.2}
\begin{tabular}{c | c | c }
$M$ & Riemannian & Complex \\
\hline
Laplacian & $\Delta_d = d\, d^{\ast} + d^{\ast} d$  & $\Delta_{\bar{\partial}}=\bpa\, \bpa^{\ast}+ \bpa^{\ast}\bpa $ \\
\hline
$\pa M$ & smooth & strongly pseudo-convex \\
\hline
Boundary & Dirichlet (D): $d(\rho \eta)~|_{\partial M}=0$ & $\bpa-$Neumann $(\bpa N)$: $\bar{\partial}^{\ast}(\rho \eta)~|_{\partial M}=0$\\
condition & Neumann (N): $d^{\ast}(\rho \eta)~|_{\partial M}=0$ & \\
\hline
Harmonic  & $\mathcal{H}_D=\{ d\eta=0, d^{\ast}\eta=0, \eta \in D\}$ & $\mathcal{H}_C=\{ \bar{\partial}\eta=0, \bar{\partial}^{\ast}\eta=0, \eta \in \bar{\partial}N\}$\\
fields& $\mathcal{H}_N=\{d\eta=0, d^{\ast}\eta=0, \eta \in N\}$ & \\
\hline
\end{tabular}
\end{table}


Clearly, the choice of suitable boundary conditions is the key to establishing Hodge theory on compact manifolds with boundary.  Hence, we ask the following question for the above symplectic Laplacians:  what type of boundary conditions on differential forms and condition on the boundary are required for a symplectic Hodge theory?

It turns out in the symplectic case, no additional condition is required of the boundary manifold besides smoothness. For the conditions on forms, we introduce the following new boundary conditions: $\partial_{+}$-Dirichlet, $\partial_{+}$-Neumann, $\partial_{-}$-Dirichlet and $\partial_{-}$-Neumann boundary conditions, denoted by $D_{+}, N_{+}, D_{-}$ and $N_{-}$, respectively.  We note in particular that $D_{+}$ and $D_{-}$ are dependent only on the symplectic structure, $\om$.   Additionally, we also introduce two boundary conditions $J-$Dirichlet and $J-$Neumann, denoted by $JD$, and $JN$, which has dependence on the almost complex structure $J$. The definitions of these boundary conditions are listed in Table 2, where $\rho$ again denotes a boundary defining function, and $d^{\Lambda}$ is the well-studied symplectic adjoint operator (see Section 2.1 for its definition). 

\begin{table}[!t]
\caption{Symplectic boundary conditions}
\renewcommand{\arraystretch}{1.2}
\begin{tabular}{c|c|c}
Notation&Definition&Relationships\\
\hline

$D_{+}$& $\partial_{+}(\rho \eta)|_{\partial M}=0$& $D \Rightarrow D_{+} $ on $\Omega^k$\\
\hline
$N_{+}$& $\partial_{+}^{\ast}(\rho \eta)|_{\partial M}=0$& $N \Leftrightarrow N_{+} $ on $P^k$ \\
\hline
$D_{-}$& $\partial_{-}(\rho \eta)|_{\partial M}=0$& $D \Rightarrow D_{-} $ on $P^k$\\
&& $D \Leftrightarrow D_{-} $ on $P^n$\\
\hline
$N_{-}$& $\partial_{-}^{\ast}(\rho \eta)|_{\partial M}=0$& $JD \Rightarrow N_{-} $ on $\Omega^k$\\
\hline
$JD$& $d^{\Lambda \ast}(\rho \eta)|_{\partial M}=0$& $JD \Rightarrow N $ on $P^k$\\
&& $JD \Leftrightarrow N_{+} $ on $P^n$\\
\hline
$JN$& $d^{\Lambda}(\rho \eta)|_{\partial M}=0$& $D \Rightarrow JN $ on $P^k$\\
&& $JN \Leftrightarrow D_{-}$ on $P^k$\\
\hline
\end{tabular}
\end{table}

The symplectic boundary conditions in Table 2 can be considered natural as they arise in the Green's formula for the corresponding Laplacians.  As is well-known, the Dirichlet boundary condition and the Neumann boundary condition appear naturally in the Green's formula of $\Delta_d$.   In the same way, for example, $D_{+}$ and $N_{+}$ arise naturally for that of $\Delta_{+}$, and $D_{-}$ and $N_{-}$ for that of $\Delta_{-}$.  Furthermore, it is worthy to note that $D_{+}$ is preserved by the differential operator $\partial_{+}$, i.e. if a form $\eta$ satisfies the $D_{+}$ boundary condition, then so does $\dpp \eta$.  Similarly, $\partial_{-}$ preserves the $D_{-}$ condition. 

These symplectic boundary conditions are also closely related to the geometry of the symplectic manifold and its boundary. For instance, we observe that the $J-$Dirichlet (Neumann) boundary condition corresponds simply to the standard Dirichlet (Neumann) boundary condition in the direction of the Reeb vector field when the boundary is of contact type.  For the $D_{+}$ condition on primitive forms, it can be thought of as the result of projecting a form with the standard Dirichlet boundary condition to its primitive component, i.e. 
\[
\pi:  \Omega^k_{D} \rightarrow P^k_{D_{+}}, k<n.
\] In other words, when a form satisfies the Dirichlet boundary conditions, its primitive component satisfies the $D_{+}$ condition. 

\begin{table}[!t]
\caption{Hodge theory on symplectic manifolds with boundary for $\Delta_{+}, \Delta_{-}, \Delta_{++}$ and $\Delta_{--}$. }
\renewcommand{\arraystretch}{1.2}
\begin{tabular}{c | c }
\hline
Laplacian & $\Delta_{+}=\dpp^{\ast}\dpp+\dpp\dpp^{\ast}$  \\
\hline
Harmonic field &  $P\mathcal{H}_{+}^k=\{ \partial_{+}\eta=0, \partial_{+}^{\ast}\eta=0\}$\\
\hline
Finite subspaces &$P\mathcal{H}_{+, D_{+}}^k, P\mathcal{H}_{+, N_{+}}^k$\\
\hline
 & $P^k=P\mathcal{H}_{+,D_{+}}^k\oplus \dpp P_{D_{+}}^{k-1}\oplus \partial_{+}^{\ast}P^{k+1}$ \\
\cline{2-2}
Decompositions &$P^k=P\mathcal{H}_{+,N_{+}}^k\oplus \dpp P^{k-1}\oplus \partial_{+}^{\ast} P_{N_{+}}^{k+1}$\\
\cline{2-2}
&$P^k=P\mathcal{H}_{+}^k \oplus \dpp P^{k-1}_{D_{+}} \oplus \partial_{+}^{\ast}  P^{k+1}_{N_{+}}$\\\hline
\end{tabular}

\vspace{10pt}

\begin{tabular}{c | c  }
\hline
Laplacian &  $\Delta_{-}=\dpm^{\ast}\dpm+\dpm\dpm^{\ast}$ \\
\hline
Harmonic field & $P\mathcal{H}_{-}^k=\{ \partial_{-}\eta=0, \partial_{-}^{\ast}\eta=0\}$\\
\hline
Finite subspaces& $P\mathcal{H}_{-, D_{-}}^k, P\mathcal{H}_{-, N_{-}}^k$\\
\hline
 & $P^k=P\mathcal{H}_{-,D_{-}}^k\oplus \partial_{-} P^{k+1}_{D_{-}}\oplus \partial_{-}^{\ast}P^{k-1}$\\
\cline{2-2}
Decompositions &$P^k=P\mathcal{H}_{-,N_{-}}^k\oplus \partial_{-} P^{k+1}\oplus \partial_{-}^{\ast}P_{N_{-}}^{k-1}$\\
\cline{2-2}
&$ P^k=  P\mathcal{H}_{-}^k \oplus \partial_{-}  P^{k+1}_{D_{-}} \oplus \partial_{-}^{\ast} P_{N_{-}}^{k-1}$\\ 
\hline
\end{tabular}


\vspace{10pt}

\begin{tabular}{c | c }
\hline
Laplacian& $\Delta_{++}=(\dpp\dpm)^{\ast}\dpp\dpm+(\dpp\dpp^{\ast})^2$  \\
\hline
Harmonic field&  $P\mathcal{H}_{++}^n=\{\eta \in P^n ~|~ \partial_{+}\partial_{-}\eta=0, \partial_{+}^{\ast}\eta=0\}$ \\
\hline
Boundary & $D_{+-}:\eta \in D_{-}, \partial_{-}\eta \in D_{+}$ \\
Condition & $N_{+-} : \eta \in N_{+}, \partial_{+}^{\ast}\eta \in N_{-}$\\
\hline
Finite subspaces &$P\mathcal{H}_{++, N_{+}}^n, P\mathcal{H}_{++, D_{+-}}^n$ \\
\hline

 & $P^n=P\mathcal{H}_{++, N_{+}}^n\oplus \partial_{+}P^{n-1}\oplus \partial_{-}^{\ast}\partial_{+}^{\ast} P_{N_{+-}}^{n}$ \\
\cline{2-2}
Decompositions & $P^n=P\mathcal{H}_{++, D_{+-}}^n\oplus \partial_{+}P^{n-1}_{D_{+}}\oplus \partial_{-}^{\ast}\partial_{+}^{\ast}P^{n}$ \\
\cline{2-2}
&$ P^n= P\mathcal{H}_{++}^n\oplus \partial_{+}P^{n-1}_{D_{+}}\oplus \partial_{-}^{\ast}\partial_{+}^{\ast} P_{N_{+-}}^{n}$ \\ 
\hline
\end{tabular}

\vspace{10pt}

\begin{tabular}{c |  c}
\hline
Laplacian & $\Delta_{--}=\dpp\dpm(\dpp\dpm)^{\ast}+(\dpm^{\ast}\dpm)^2$ \\
\hline
Harmonic field & $P\mathcal{H}_{--}^n=\{\eta \in P^n ~|~ \partial_{-}\eta=0, \partial_{-}^{\ast}\partial_{+}^{\ast}\eta=0\}$\\
\hline

Finite subspaces  & $P\mathcal{H}_{-, D_{-}}^n, P\mathcal{H}_{--, N_{+-}}^n$\\
\hline
 & $P^n=P\mathcal{H}_{--, D_{-}}^n\oplus \partial_{+}\partial_{-} P_{D_{+-}}^{n}\oplus \partial_{-}^{\ast}P^{n-1}$\\
\cline{2-2}
Decompositions  &$P^n=P\mathcal{H}_{--, N_{+-}}^n\oplus \partial_{+}\partial_{-}  P^n\oplus \partial_{-}^{\ast}P_{N_{-}}^{n-1}$\\
\cline{2-2}
 &$P^n=P\mathcal{H}_{--}^n\oplus \partial_{+}\partial_{-} P_{D_{+-}}^{n}\oplus \partial_{-}^{\ast}P_{N_{-}}^{n-1}$\\ 
\hline
\end{tabular}
\end{table}

To obtain the Hodge theory with the above symplectic boundary conditions, we make use of the theory of elliptic boundary value problems (BVPs) \cite{H3}.  The boundary conditions involved will then be elliptic in the sense of Lopatinski-Shapiro. We generalize the argument in \cite{SM} to prove the smoothness of weak solutions under certain assumptions and then apply the theory of elliptic boundary value problems to obtain the desired Hodge theory.  We note that the necessity to consider the fourth-order symplectic Laplacians, i.e. $\Delta_{++}$ and $\Delta_{--}$, is rather special for the symplectic case.  In Table 3, we list our main results for the Hodge theory of symplectic Laplacians.

The Hodge theory established in this paper have various applications.  For one, we use it to identify isomorphisms between symplectic cohomologies and spaces of harmonic fields of symplectic Laplacians with certain boundary conditions.  These isomorphisms show that certain symplectic cohomologies are indeed still finite-dimensional when the boundary is not vanishing.  Furthermore, the dimensions of certain harmonic fields with boundary conditions can be regarded as invariants of the symplectic structure.  In a different direction, we also utilize Hodge theory to obtain various Poincar\'{e} lemmas and solve a number of boundary value problems related to the existence of harmonic fields. As a consequence,  we show that the spaces of symplectic harmonic fields with no boundary conditions are infinite-dimensional.


The structure of this paper is as follows. In Section 2, we provide some basic definitions and lemmas needed in this paper. In Section 3, we define the symplectic boundary conditions and discuss some of their properties. In Section 4, we obtain the finiteness and decompositions results about the Hodge theory for symplectic Laplacians through the theory of elliptic boundary value problems. In Section 5, we study the symplectic cohomologies of compact symplectic manifolds with boundary. We build isomorphisms between these cohomologies and subspaces of harmonic fields with certain boundary conditions. In Section 7, we apply our results to boundary value problems and prove various  Poincar\'{e} lemmas, which give necessary and sufficient conditions for a form to be $\dpp$, $\dpm$, or $\dpp\dpm$-exact on a compact symplectic manifold with boundary. We also prove the infiniteness of the spaces of harmonic fields without boundary conditions. In Section 7, we discuss a few additional observations of symplectic boundary conditions which may point further to other invariants on symplectic manifolds with boundary.

Let us note that a standard method to establish the Hodge theory is through the so-called Gaffney inequality, or more often called the G\"{a}rding's inequality in the case of manifolds without boundary.   In our study, we initially took the Gaffney inequality approach to symplectic Laplacians which helped us gain intuition concerning the Laplacians and boundary conditions.  However, the results obtained from this approach thus far involved much stronger and more complicated boundary conditions than those appearing in Table 3.  But with its relevance for analysis, we deem it still worthwhile to include some of our Gaffney inequalities results for the second-order symplectic Laplacians as part of an Appendix.

\

\noindent{\it Acknowledgements.~} 
We would like to thank X. Dai, R.-T. Huang, L. Ni, Y. S. Poon, M. Schecter, C.-J. Tsai, G. Xu, and S.-T. Yau for helpful comments and discussions.   Additionally, we are grateful to S.-Y. Li, Z. Lu, C.-L. Terng and especially P. Li for their interest, encouragement, and input in this work.

\section{Preliminaries }
We recall some basic definitions and properties in symplectic geometry and Riemannian geometry, cf.\cite{TY1, TY2}.  Lemmas and propositions given there will be stated here without proof. 
\subsection{Primitive structures on symplectic manifolds}
Given a symplectic manifold $(M^{2n}, \omega)$,  let $\Omega^k$ denote the space of smooth $k-$forms on $M$. With respect to local coordinates, write the symplectic form as $\omega= \frac{1}{2}\sum \omega_{ij}~ dx^{i}\wedge dx^{j}$. The Lefschetz operator $L$ and its dual operator $\Lambda$ are defined by
\begin{align*}
L: \Omega^k \rightarrow \Omega^{k+2}, L(\eta)&=\omega \wedge \eta,\\
\Lambda: \Omega^k \rightarrow \Omega^{k-2}, \Lambda(\eta)&=\frac{1}{2}(\omega^{-1})^{ij}i_{\partial_{x^i}}i_{\partial_{x^j}}\eta
\end{align*} where $i$ denotes the interior product, and $\omega^{-1}$ is the inverse matrix of $\omega$. Define the degree counting operator 
\begin{equation*}
H=\underset{k}{\sum}(n-k) \prod ^k
\end{equation*} with $\prod^k: \Omega^{\ast} \rightarrow \Omega^k$ as the projection operator onto forms of degree $k$. As is known, $L$ and $\Lambda$ together with $H$ give a representation of $sl(2)$ algebra acting on $\Omega^{\ast}$:
\begin{equation*}
[\Lambda, L]=H, \, [H, \Lambda]=2\Lambda, \, [H,L]=-2L.
\end{equation*}

This $sl(2)$ representation allows a "Lefschetz" decomposition of forms in terms of irreducible finite-dimensional $sl(2)$ modules. The highest weight states of these irreducible $sl(2)$ modules are the spaces of primitive forms, denoted by $P^{\ast}$.
\begin{definition}
A form $\eta \in \Omega^k$ is called primitive if $\Lambda \eta=0$. This is equivalent to the condition $L^{n-k+1}\eta=0$. 
\end{definition} 
As implied by the definition, the degree of the primitive form is constrained to be $k \leq n$. Given $\eta \in \Omega^k$, there is a unique Lefschetz decomposition into primitive forms as 
\[
\eta =\underset{r \geq \max(k-n, 0)}{\sum}\frac{1}{r!}L^rB_{k-2r}.
\] Here each $B_{k-2r}\in P^{k-2r}$ can be expressed in terms of $\eta$: $B_{k-2r}=\left(\underset{s=0}{\sum}a_{r,s}\frac{1}{s!}L^s\Lambda^{r+s}\right)\eta$. Thus, each term of this decomposition can be labeled by a pair $(r,s)$ corresponding to the space 
 \[
 \mathcal{L}^{r,s}=\left\{A\in \Omega^{2r+s}~|~ A=L^rB_s \,\text{with}\, B_s \in P^s\right\}.
 \] We cite the following lemma about $L, \Lambda$ and $H$ from \cite{TY2}. 
 \begin{lemma}
On symplectic  manifolds, the following relations hold:
\begin{itemize}
\item $[\Lambda, L^r]=(H+r-1)rL^{r-1}$ for $r\geq 1$;
\item $L\Lambda=(H+R+1)R$;
\item $\Lambda L=(H+R)(R+1)$.
\end{itemize} Here, the operator $R$ is defined as $R(L^rB_s)=rL^rB_s$.
\end{lemma}

\subsection{Differential operators $\partial_{+}, \partial_{-}$, and $d^{\Lambda}$}
 We consider the action of the exterior derivative operator $d$ on $\mathcal{L}^{r,s}$, cf \cite{TY2}. 
\begin{prop}
$d$ acting on $\mathcal{L}^{r,s}$ leads to at most two terms:
\begin{equation*}
d: \mathcal{L}^{r,s} \rightarrow \mathcal{L}^{r,s+1}\oplus \mathcal{L}^{r+1,s-1}
\end{equation*} with 
\begin{align*}
dL^rB_s&=L^r(dB_s)=L^rB_{s+1}+L^{r+1}B_{s-1}\, \text{ when $s<n$},\\
dL^rB_n&=L^r(dB_n)=L^{r+1}B_{n-1}.
\end{align*} 
\end{prop} This result is implied by the closeness of the symplectic form $\omega$,  and the following formulas:
\begin{itemize}
  \item If $s <n$, $dB_s=B_{s+1}+LB_{s-1}$,
  \item If $s=n$, $dB_n=LB_{n-1}$.
\end{itemize} By this proposition, we can define the decomposition of $d$ into two linear differential operators $(\partial_{+}, \partial_{-})$. 
\begin{definition}
On a symplectic manifold $(M, \omega)$, we define the first order differential operators $\partial_{+}, \partial_{-}$ by the property:
\begin{align*}
&\partial_{+}:\mathcal{L}^{r,s} \rightarrow \mathcal{L}^{r,s+1}, \qquad\partial_{+}(L^{r}B_{s})=L^{r}B_{s+1},\\
&\partial_{-}:  \mathcal{L}^{r,s} \rightarrow \mathcal{L}^{r,s-1}, \qquad \partial_{-}(L^{r}B_{s})=L^{r}B_{s-1}
\end{align*}  such that $d=\partial_{+}+L\partial_{-}$. Here $B_s, B_{s+1}, B_{s-1} \in P^{\ast}$ and $dB_s=B_{s+1}+LB_{s-1}$.
\end{definition} 
When acting on primitive forms, $\partial_{+}$ and $\partial_{-}$ can be equivalently written as follows: 
\begin{lemma}
Acting on primitive differential forms, operators $\partial_{+}, \partial_{-}$ have the following expressions:
\begin{align*}
\partial_{+}&=d-LH^{-1}\Lambda d,\\
\partial_{-}&=H^{-1}\Lambda d.
\end{align*}
\end{lemma}

In explicit calculations (e.g. in Appendix A) , it can be useful to modify the differential operator $\partial_{-}$ to reduce the number of constant factors that arise.  We define 
\[
 \dpmm=(H+R)\partial_{-}.
\] Thus, 
\[
d=\partial_{+}+L(H+R)^{-1}\dpmm.
\] With this definition, the following properties hold:
\begin{prop}
On $(M^{2n}, \omega)$, the symplectic differential operators $(\partial_{+}, \dpmm)$ satisfy:
\begin{itemize}
\item $\partial_{+}^2=(\dpmm)^2=0$,
\item $(H+R)\partial_{+}\dpmm=(H+R+1)\dpmm\partial_{+}$ on $\mathcal{L}^{r,s}$,
\item $[\partial_{+}, L]=[L\,\partial_{-}, L]=0$. But $[L\,\dpmm, L]=-L^2\partial_{-} $.
\end{itemize}
\end{prop}
Besides $d$, there is another first-order differential operator of interest in this paper $$d^{\Lambda}=d\,\Lambda - \Lambda\,d : \Omega^{k} \rightarrow \Omega^{k-1}.$$ 
With $d$ and $d^{\Lambda}$, $\partial_{+}$ and $\partial_{-}$ can be expressed as follows.
\begin{lemma}\label{ex}
On a symplectic manifold $(M, \omega)$, $\partial_{+}$ and $\partial_{-} $ can be expressed as 
\begin{align*}
\partial_{+}&=\frac{1}{H+2R+1}\left[ (H+R+1)d+Ld^{\Lambda}\right],\\
\partial_{-}&=\frac{1}{(H+2R+1)(H+R)}\left[\Lambda d-(H+R)d^{\Lambda}\right].
\end{align*}
\end{lemma}

Let us also note the following proposition from \cite{TY1}.
\begin{prop}
With respect to the $sl(2)$ representation $(L, \Lambda, H)$, the differential operators $(d, d^{\Lambda}, dd^{\Lambda})$ satisfy the following commutation relations:
\begin{align*}
&[d, L]=0, && [d, \Lambda]=d^{\Lambda}, &&& [d, H]=d,\\
&[d^{\Lambda}, L]=d, && [d^{\Lambda}, \Lambda]=0, &&&[d^{\Lambda}, H]=-d^{\Lambda},\\
&[dd^{\Lambda}, L]=0,&&[dd^{\Lambda}, \Lambda]=0,&&&[dd^{\Lambda}, H]=0.
\end{align*}
\end{prop}

\subsection{Conjugate relations}
Let $(\omega, J, g)$ be a compatible triple on the symplectic manifold $(M, \omega)$ with $J$ as an almost complex structure and $g$ as a Riemannian metric. With respect to the almost complex structure $J$, there is the decomposition $\Omega^k=\underset{p+q=k}{\oplus}\Omega^{p,q}$. Then define the operator 
\[
\mathcal{J}=\underset{p,q}{\sum}(\sqrt{-1})^{p-q}\prod^{p,q}
\] which projects a $k-$form onto its $(p,q)$ parts timing with the multiplicative factors $(\sqrt{-1})^{p-q}$. Notice $\mathcal{J}^2=(-1)^k$ acting on $k-$forms. The operator $\mathcal{J}$ communicates with $L$ and $\Lambda$.
\begin{lemma}\label{jl}
 For a triple $(\omega, \mathcal{J}, g)$, there is
\begin{equation*}
[\mathcal{J}, L]=0,\,[\mathcal{J}, \Lambda]=0.
\end{equation*}
\end{lemma}
This is because that the symplectic form $\omega$ is a $(1,1)$-form with respect to the almost complex structure $J$. Moreover, the operator $\mathcal{J}$ defines the following conjugate relations (\cite{TY1,TY2}) between differential operators.
\begin{lemma}\label{conjugate}
For a compatible triple $(\omega, \mathcal{J}, g)$ on a symplectic manifold, let $d^{\ast}, d^{\Lambda \ast}, \partial_{+}^{\ast}$ and $\dpmm^{\ast}$ be the adjoint operators of the corresponding differential operators, respectively. Then there are conjugate relations:
\begin{itemize}
\item $d^{\Lambda}=-\mathcal{J}d^{\ast}\mathcal{J}^{-1}$ and $d^{\Lambda \ast}=-\mathcal{J}d\mathcal{J}^{-1}$;
\item $\dpmm=\mathcal{J}\partial_{+}^{\ast}\mathcal{J}^{-1}$ and $\dpmm^{\ast}=\mathcal{J}\partial_{+}\mathcal{J}^{-1}$.
\end{itemize} 
\end{lemma} 
Define the $d^{\Lambda}$ Laplacian: $\Delta_{d^{\Lambda}}=d^{\Lambda\ast}d^{\Lambda}+d^{\Lambda}d^{\Lambda\ast}$.   Then Lemma \ref{conjugate} implies that this operator is conjugate to the Laplace operator $\Delta_d$.  
\begin{cor} Let $(\omega, \mathcal{J}, g)$ be a compatible triple on a symplectic manifold. Then the following conjugate relation holds:
\begin{equation*}
\Delta_{d^{\Lambda}}=\mathcal{J}\,\Delta_d\,\mathcal{J}^{-1}.
\end{equation*} 
\end{cor}
Because of this conjugate relation, the ellipticity of $\Delta$ implies that  of $\Delta_{d^{\Lambda}}$. Moreover, we have the following expressions of adjoint operators according to Lemma \ref{ex}.
\begin{lemma}\label{c2}
On a symplectic manifold $(M, \omega)$ with a compatible Riemannian metric $g$, the adjoints $(\partial^{\ast}_{+}, \partial^{\ast}_{-})$ take the forms
\begin{align*}
\partial^{\ast}_{+}&=[d^{\ast}(H+R+1)+d^{\Lambda\ast}\Lambda](H+2R+1)^{-1},\\
\partial^{\ast}_{-}&=[d^{\ast}(H+R+1)^{-1}L-d^{\Lambda\ast}](H+2R+1)^{-1}.
\end{align*}
\end{lemma}

\begin{cor}\label{c1}
On $P^k$, the adjoints $(\partial^{\ast}_{+}, \partial^{\ast}_{-})$ take the forms
\begin{align*}
\partial^{\ast}_{+}&=d^{\ast},\\
\partial^{\ast}_{-}&=(n-k)^{-1}d^{\ast}L-(n-k+1)^{-1}Ld^{\ast}.
\end{align*} Moreover,
\begin{align*}
\partial^{'\ast}_{-}&=d^{\ast}L-\frac{n-k}{n-k+1}Ld^{\ast}.
\end{align*}

\end{cor}
\subsection{Symplectic Laplacians}
On a symplectic manifold $(M, \omega)$, there exists an elliptic complex on primitive spaces  \cite{TY2} (see also \cite{Smith, Eastwood, E1}):
$$\begin{CD}
0@>\partial_{+}>>P^0@>\partial_{+}>>P^1@>\partial_{+}>>\cdots @>\partial_{+}>>P^{n-1}@>\partial_{+}>>P^n\\
    @. @. @. @.   @.                                                                                                            @VV{\partial_{+}\partial_{-}}V\\
0@<\partial_{-}<<P^1@<\partial_{-}<<P^2@<\partial_{-}<<\cdots @<\partial_{-}<<P^{n-1}@<\partial_{-}<<P^n
\end{CD}$$
Note that a special part of this complex is the second-order differential operators $\dpp\dpm$ acting on the middle degree primitive forms, $P^n$. We define the following Laplacians  associated to this elliptic complex:
\begin{align*}
\Delta_{+}&=\partial_{+}\partial_{+}^{\ast}+\partial_{+}^{\ast}\partial_{+}, \, \text{on $P^k$, for $k<n$};\\
\Delta_{-}&=\partial_{-}\partial_{-}^{\ast}+\partial_{-}^{\ast}\partial_{-},\, \text{on $P^k$, for $k<n$};\\
\Delta_{++}&=(\partial_{+}\partial_{-})^{\ast}(\partial_{+}\partial_{-})+(\partial_{+}\partial_{+}^{\ast})^2, \, \text{on $P^n$};\\
\Delta_{--}&=(\partial_{+}\partial_{-})(\partial_{+}\partial_{-})^{\ast}+(\partial_{-}^{\ast}\partial_{-})^2, \, \text{on $P^n$}.
\end{align*}For calculational simplification, it is sometimes useful to define the following Laplacian by replacing $\partial_{-}$ by $\partial_{-}'$ in the definition above.
\begin{align*}
\Dmm&=\dpmm\dpmm^{\ast}+\dpmm^{\ast}\dpmm,\, \text{on $P^k$, for $k<n$}.
\end{align*} 
The ellipticity of the complex implies that the operators $\Delta_{+}$ and $\Delta_{-}$ are elliptic on $P^k$ for $k<n$.  Similarly $\Delta_{-}'$ is elliptic on $P^k$ when $k<n$. The ellipticity of the Laplacians $\Delta_{++}$ and $\Delta_{--}$ however may not be immediately obvious.  We will prove that $\Delta_{++}$ and $\Delta_{--}$ are elliptic on $P^n$ explicitly by calculating their symbols.  But do so, it is useful to utilize two other fourth-order symplectic Laplacians on $\Omega^k$ that were introduced in \cite{TY1} and modified here as follows.
\begin{definition}
For any $\eta \in \Omega^k$, define the following operators:
\begin{align*}
\Delta_{dd^{\Lambda}}(\eta)&=d^{\Lambda\ast}d^{\ast}dd^{\Lambda} \eta+\frac{1}{4}\left(dd^{\ast}+d^{\Lambda}d^{\Lambda \ast}\right)^2\\
\Delta_{d+d^{\Lambda}}(\eta)&=dd^{\Lambda}d^{\Lambda\ast}d^{\ast} \eta+\frac{1}{4}\left(d^{\ast}d+d^{\Lambda \ast}d^{\Lambda}\right)^2.
\end{align*}
\end{definition}
\begin{prop}\label{sb}
The operators $\Delta_{dd^{\Lambda}}$ and $\Delta_{d+d^{\Lambda}}$ are elliptic. Moreover, they reduce to $\Delta_{++}$ and $\Delta_{--}$ on $P^n$, respectively. Therefore, $\Delta_{++}$ and $\Delta_{--}$ are elliptic on $P^n$.
\end{prop}
\begin{proof}
Fix a point $x \in M$ and let $\xi \in \Omega^1_x$ be any normalized $1-$form. Choose a basis $\{w_i\}$ of the cotangent space at $x$ such that $w_1=\xi$ and symplectic form takes the form $\omega= w_1\wedge w_2+\cdots+w_{2n-1} \wedge w_{2n}$. Let $\{e_i\}$ denote the dual basis. Then any $k$-form $\eta$ can be expressed in the form:
\[
\eta= w_1\wedge \beta_1+w_2\wedge \beta_2+ w_{12}\wedge \beta_3+\beta_4.
\] Here $\beta_i$ are forms without neither $w_1$ nor $w_2$ in their components, and $w_{ij}$ denotes $w_i \wedge w_j$. We have the following symbol calculations at $x$:
\begin{align*}
\sigma(d)(\xi)\eta&=w_1\wedge \eta=w_{12}\wedge\beta_2+w_1\wedge \beta_4,\\
\sigma(d^{\ast})(\xi)\eta&=-i_{e_1}\eta=-(\beta_1+w_2\wedge\beta_3),\\
\sigma(d^{\Lambda})(\xi)\eta&=w_1\wedge \beta_3-\beta_2=-i_{e_2}\eta,\\
\sigma(d^{\Lambda \ast})(\xi)\eta&=-w_{12}\wedge\beta_1+w_2\wedge \beta_4=w_2\wedge \eta.
\end{align*}
Here, $i_{v}$ denotes the interior product with the tangent vector $v$.
Therefore, we obtain
\begin{align*}
\sigma(\Delta_{dd^{\Lambda}})(\xi)\eta&=w_2\wedge \beta_2+w_1\wedge \beta_1+\frac{1}{4}w_{12}\wedge \beta_3+\frac{1}{4}\beta_4;\\
\sigma(\Delta_{d+d^{\Lambda}})(\xi)\eta&=w_1\wedge \beta_1+w_2\wedge \beta_2+\frac{1}{4}w_{12}\wedge\beta_3+\frac{1}{4}\beta_4.
\end{align*} These explain the ellipticity of both $\Delta_{dd^{\Lambda}}$ and $\Delta_{d+d^{\Lambda}}$.
For the last claim,  take $\eta \in P^n$ and we get
\begin{align*}
d\eta&=L\partial_{-}\eta, \, \partial_{+}^{\ast}\eta=d^{\ast}\eta\\
d^{\Lambda}\eta&=-\partial_{-}\eta,\, d^{\Lambda \ast}\eta=Ld^{\ast}\eta=L\partial_{+}^{\ast}\eta.
\end{align*}It is not hard so see that
\begin{align*}
\partial_{+}\partial_{-}\eta&=-dd^{\Lambda}\eta\\
\partial_{+}\partial_{+}^{\ast}\eta&=\frac{1}{2}(dd^{\ast}+d^{\Lambda}d^{\Lambda \ast})\eta.
\end{align*} Thus $\Delta_{dd^{\Lambda}}=\Delta_{++}$ on $P^n$. Similarly,  $\Delta_{d+d^{\Lambda}}=\Delta_{--}$ on $P^n$.
\end{proof}
\begin{remark}
Generally $\Delta_{++}$ and $\Delta_{--}$ are not elliptic on $P^k $ when $k<n$. 
\end{remark}

\section{Symplectic boundary conditions}
 Given a compact symplectic manifold $(M^{2n}, \omega)$ with smooth boundary $\partial M$, let $(\omega, J, g)$ be a compatible triple on it. Let $\rho$ be a  boundary defining function, i.e.
\begin{itemize}
\item $\rho < 0$ on $M$ and $\rho(x)=0$ if and only if $x \in \partial M$,
\item the norm of  gradient $|\nabla \rho| =1$ on $\partial M$.  
\end{itemize} 
 \begin{definition}
For a form $\eta$ which is well defined along $\partial M$, we say $\eta$ satisfies
\begin{itemize}
\item Dirichlet boundary condition, denoted by $\eta \in D$, if $d(\rho \eta)|_{\partial M}=0$. 
\item Neumann boundary condition, denoted by  $\eta \in N$, if $d^{\ast}(\rho \eta)|_{\partial M}=0$. 
\item $J$-Dirichlet boundary condition, denoted by  $\eta \in JD$, if $d^{\Lambda \ast}(\rho \eta)|_{\partial M}=0$.
\item $J$-Neumann boundary condition, denoted by  $\eta \in JN$, if $d^{\Lambda}(\rho \eta)|_{\partial M}=0$.\item $\partial_{+}$-Dirichlet boundary condition, denoted by  $\eta \in D_{+}$, if $\partial_{+}(\rho \eta|_{\partial M})=0$.
\item $\partial_{+}$-Neumann boundary condition, denoted by  $\eta \in N_{+}$, if $\partial_{+}^{\ast}(\rho \eta|_{\partial M})=0$. 
\item $\partial_{-}$-Dirichlet boundary condition, denoted by  $\eta \in  D_{-}$, if $\partial_{-}(\rho \eta|_{\partial M})=0$.
\item $\partial_{-}$-Neumann boundary condition, denoted by  $\eta \in  N_{-}$, if if $\partial_{-}^{\ast}(\rho \eta|_{\partial M})=0$. 
\end{itemize}
\end{definition} The $J$-Dirichlet and $J$-Neumann boundary condition are named based on the following lemma.
\begin{lemma}
For a form $\eta \in \Omega^k$, there are
\begin{itemize}
\item  $\eta \in JD$ if and only if $\mathcal{J}\eta \in D$,
\item  $\eta \in JN$ if and only if $\mathcal{J}\eta \in N$.
\end{itemize}
\end{lemma}
\begin{proof}
Since $d^{\Lambda \ast}=\mathcal{J}d\mathcal{J}^{-1}$, there is 
\[
\eta \in JD \, \text{if and only if }\, d^{\Lambda \ast}(\rho \eta)|_{\partial M}=0\, \text{if and only if }\, d(\rho \mathcal{J}\eta)|_{\partial M}=0.
\] This means that $\eta \in JD$ if and only if $\mathcal{J}\eta \in D$. Similarly, the relation $d^{\Lambda}=\mathcal{J}d^{\ast}\mathcal{J}^{-1}$ implies that $\eta \in JN$ if and only if $\mathcal{J}\eta \in N$.
\end{proof}

In fact, these boundary conditions are natural in the sense that they appear in various Green's formulas.  We recall the following property \cite{T}:
\begin{lemma}[Green's formula for first-order differential operators]
If $M$ is a smooth, compact manifold with boundary and $P$ is a first-order differential operator acting on sections of a vector bundle, then
\begin{equation*}
(Pu,v)-(u,P^t v)=\int_{\partial M}\langle \sigma_p(x,\overrightarrow{\!\!n})u, v \rangle
\end{equation*} with $P^t$ as the dual operator of $P$, $\sigma_p$ as the symbol of $P$ and $\overrightarrow{\!\!n}$ as the outward normal along the boundary $\partial M$.
\end{lemma}

For example, for the exterior differential operator $d$ acting on $\Omega^{\ast}(M)$, 
\begin{align*}
(d\eta,w)-(\eta, d^{\ast}w)&=\int_{\partial M}\langle \sigma_{d}(x, \overrightarrow{\!\!n}) \eta, w\rangle\\
(d\eta,w)-(\eta, d^{\ast}w)&=-\int_{\partial M}\langle \eta, \sigma_{d^{\ast}}(c,\overrightarrow{\!\!n})w\rangle
\end{align*}
are implied by the proposition. Here
\[
\sigma_d(x,\overrightarrow{\!\!n}) \eta=d(\rho \eta)\,\text{and}\, \sigma_{d^{\ast}}(x,\overrightarrow{\!\!n}) w=d^{\ast}(\rho w).
\] This is the standard result of the Green's formula for $d$ and $d^{\ast}$.  Additionally,

\begin{cor}[Green's formula for $d^{\Lambda}$, $\partial_{+}$ and $\partial_{-}$]\label{green}
Given a compact symplectic manifold $M$ with smooth boundary $\partial M$, let $(\omega, J, g)$ be a compatible triple on it. Then for any $\phi, \psi \in \Omega^k$, there are
\begin{align*}
(d^{\Lambda}\phi, \psi)-(\phi, d^{\Lambda \ast}\psi)&=\int_{\partial M}\langle d^{\Lambda}(\rho \phi), \psi\rangle=-\int_{\partial M}\langle \phi, d^{\Lambda \ast}(\rho \psi)\\
(\partial_{+}\phi, \psi)-(\phi, \partial^{\ast}_{+}\psi)&=\int_{\partial M}\langle \partial_{+}(\rho \phi), \psi\rangle=-\int_{\partial M}\langle \phi, \partial_{+}^{\ast}(\rho \psi)\\
(\partial_{-}\phi, \psi)-(\phi, \partial_{-}^{\ast}\psi)&=\int_{\partial M}\langle \partial_{-}(\rho \phi), \psi \rangle=-\int_{\partial M}\langle \phi, \partial_{-}^{\ast}(\rho \psi)\rangle.
\end{align*}
\end{cor} 
\begin{proof}
This is a direct application of the general Green's formula. We only need to point out that 
\begin{align*}
\sigma_{d^{\Lambda}}(x, \overrightarrow{\!\!n}) \phi&=d^{\Lambda}(\rho \phi)\\
\sigma_{\partial_{+}}(x, \overrightarrow{\!\!n}) \phi&=\partial_{+}(\rho \phi)\\
\sigma_{\partial_{-}}(x, \overrightarrow{\!\!n}) \phi&=\partial_{-}(\rho \phi)
\end{align*}on $\partial M$ by definition.
\end{proof}

\begin{remark}
Notice $\partial_{-}'(\rho \eta)|_{\partial M}=0$ is equivalent to $\partial_{-}(\rho \eta)|_{\partial M}=0$, and $\partial_{-}^{'\ast}(\rho \eta)|_{\partial M}=0$ is equivalent to $\partial_{-}^{\ast}(\rho \eta)|_{\partial M}=0$.\end{remark}
From the definition, the following adjoint relations between the different symplectic boundary conditions hold.
\begin{lemma}
Let $\ast$ be the Hodge star operator and $\eta \in \Omega^k$. Then 
\begin{align*}
&\eta \in D \, \text{if and only if }\, \ast\eta \in N\\
&\eta \in JD \, \text{if and only if }\, \ast\eta \in JN\\
&\eta \in D_{+} \, \text{if and only if }\, \ast\eta \in N_{+}\\
&\eta \in D_{-} \, \text{if and only if }\, \ast\eta \in N_{-}.
\end{align*}
\end{lemma}
Now we are ready to explain the relations in Table 2 of the Introduction with following lemma.
\begin{lemma}\label{bc}
For any $\eta \in \Omega^k$, 
\[
\eta \in D\Rightarrow \eta \in D_{+}, \eta \in D_{-}, \qquad \eta \in JD \Rightarrow \eta \in N_{+}, \eta \in N_{-}.
\] When $\eta \in P^k$, 
\[
 \eta \in N \Leftrightarrow \eta \in N_{+}, \qquad \eta \in JN \Leftrightarrow \eta \in D_{-}. 
\] Moreover,
\begin{align*}
\eta &\in D \Rightarrow \eta \in JN, \quad \eta \in JD \Rightarrow \eta \in N\quad \text{ when $\eta \in P^k$}\\
\eta &\in D\Leftrightarrow \eta \in JN, \quad \eta \in JD \Leftrightarrow \eta \in N \quad \text{ when $\eta \in P^n$}.
\end{align*}
\end{lemma} 
\begin{proof}
Since $d=\partial_{+}+L\partial_{-}$, it is easy to see that $\eta \in D$ implies $\eta \in D_{+}$ and $\eta \in D_{-}$. If $\eta \in JD$, then $\mathcal{J}\eta \in D$ which implies that $\mathcal{J} \eta \in D_{+}$ and $\mathcal{J}\eta \in D_{-}$. By the conjugate relations, we have 
\begin{align*}
\partial_{+}(\rho \mathcal{J}\eta)|_{\partial M}=0 & \Rightarrow \mathcal{J}\partial_{+}(\mathcal{J}\rho\eta)|_{\partial M}=0 \Rightarrow \partial_{-}^{\ast}(\rho\eta)|_{\partial M}=0 \\
\partial_{-}(\rho \mathcal{J}\eta)|_{\partial M}=0 & \Rightarrow \mathcal{J}\partial_{-}(\rho \mathcal{J}\rho \eta)|_{\partial M}=0 \Rightarrow \partial_{+}^{\ast}(\rho\eta)|_{\partial M}=0.
\end{align*} Therefore, $\eta \in JD$ implies both $\eta \in N_{+}$ and $\eta \in N_{-}$. 

Since $d^{\ast}=\partial_{+}^{\ast}$ and $d^{\Lambda}=-H^{-1}\partial_{-}$ on $P^k$, the equivalences $\eta \in N \Leftrightarrow \eta \in N_{+}, \, \eta \in JN \Leftrightarrow \eta \in D_{-}$ hold for $\eta \in P^k$. 
The results in the case of $\eta \in P^n$ are immediate.
\end{proof}

Next we will illustrate that these boundary conditions can be preserved by certain differential operators. Fix a point $x \in \partial M$ and choose a normal basis $\{w_i\}$ of the cotangent space at $x$, such that $w_{1}=d\rho$ and $\omega=\underset{i}{\sum}w_{2i-1}\wedge w_{2i}$. Then $w_{2}=\mathcal{J} w_{1}$ by the compatible condition. Let $\{e_i\}$ be the dual basis. We obtain the following local characterization of these boundary conditions with respect to $\{w_i\}$:
\begin{lemma}
At a point $x\in \partial M$, with respect to the basis $\{w_i\}$ as above, we have 
\begin{itemize}
\item $\eta \in D$ if and only if $w_1 \wedge \eta=0$ along $\partial M$,
\item $\eta \in N$ if and only if $i_{e_1} \eta=0$ along $\partial M$,
\item $\eta \in JD$ if and only if $w_2 \wedge \eta=0$ along $\partial M$,
\item $\eta \in JN$ if and only if $i_{e_2} \eta=0$ along $\partial M$.
\end{itemize}
\end{lemma} 
\begin{proof}
By the definitions of $\rho, d$ and $d^{\ast}$, it is not hard to see that the first two statements hold. For the other two, we have
\begin{itemize}
\item $\eta \in JD$ if and only if $d(\rho \mathcal{J}\eta)~|_{\partial M}=0$ if and only if $w_1\wedge \mathcal{J}\eta ~|_{\partial M}=0$ which is equivalent to $w_2\wedge \eta =0$ along $\partial M$;
\item $\eta \in JN$ if and only if $d^{\ast}(\rho \mathcal{J}\eta)~|_{\partial M}=0$ if and only if $i_{e_1} \mathcal{J}\eta~|_{\partial M}=0$ which is equivalent to $i_{e_2}\eta =0$ along $\partial M$.
\end{itemize}
\end{proof}
By definitions, the Dirichlet boundary condition and the Neumann boundary condition are defined with respect to the direction of the outer normal $d\rho$.  Equivalently then, the $J-$Dirichlet boundary condition and $J-$Neumann boundary condition are defined with respect to the direction given by $\mathcal{J} d\rho$. Moreover, by the compatibility of $(\omega, J, g)$, these two directions $d\rho$ and $\mathcal{J}d\rho$ are orthogonal:
\[
g(d\rho, \mathcal{J}d\rho)=\omega(d\rho, \mathcal{J}^2d\rho)=-\omega(d\rho, d\rho)=0.
\] 

Define the tangential projection and the normal projection of $\eta$ along $\partial M$ as follows:
\begin{equation*}
\eta^t=i_{e_{1}} (w_{1}\wedge \eta), \, \eta^n=w_{1}\wedge(i_{e_{1}} \eta).
\end{equation*} Then $\eta \in D$ if and only if $\eta^t=0$, and $\eta \in N$ if and only if $\eta^n=0$. We cite the following proposition from \cite{S}.
\begin{prop}With assumptions above, the following results hold.
\begin{itemize}
\item The normal and tangential components of $\eta \in \Omega^k$ are Hodge adjoint to each other:
\begin{equation*}
\ast(\eta^n)=(\ast \eta)^t, \, \ast(\eta^t)=(\ast \eta)^n.
\end{equation*} Here, $\ast( \eta^n)$ and $\ast(\eta^t)$ are understood by the action of $\ast$ on arbitrary extension of $\eta^n$ and $\eta^t$, respectively, followed by the restriction of $\partial M$.
\item The exterior derivative commutes with tangential projection, and the co-differential with normal projection of $\eta \in \Omega^k$ in the following sense:
\begin{equation*}
(d \eta)^t= d(\eta^t), \, (d^{\ast}\eta)^n=d^{\ast} (\eta^n)
\end{equation*}
\end{itemize}
\end{prop} 
This proposition tells us that the boundary condition $D$ is dual to the condition $N$ by the operator $\ast$. An immediate consequence of this proposition is as follows.
\begin{cor}
For a form $\eta \in \Omega^k$, 
\begin{align*}
\eta \in D& \Rightarrow d\eta \in D, \qquad \eta \in N\Rightarrow d^{\ast}\eta \in N;\\
\eta \in JD&\Rightarrow d^{\Lambda\ast}\eta \in JD, \qquad \eta \in JN \Rightarrow d^{\Lambda} \eta\in JN.
\end{align*}
When $\eta \in P^k$, 
\begin{align*}
\eta \in D_{+} &\Rightarrow \partial_{+}\eta \in D_{+}, \qquad \eta \in N_{+} \Rightarrow \partial_{+}^{\ast} \eta \in N_{+};\\
\eta \in D_{-} &\Rightarrow \partial_{+}\eta \in D_{-}, \qquad \eta \in N_{-} \Rightarrow \partial_{-}^{\ast} \eta \in N_{-}.
\end{align*}
\end{cor}
\begin{proof}
The first line is obvious by the proposition above. For the second line, there are the relations
\begin{align*}
\eta \in JD \Rightarrow \mathcal{J}\eta \in D\Rightarrow d\mathcal{J}\eta \in D \Rightarrow d^{\Lambda\ast}\eta \in JD,\\
\eta \in JN \Rightarrow \mathcal{J}\eta \in N\Rightarrow d\mathcal{J}\eta \in N \Rightarrow d^{\Lambda}\eta \in JN.
\end{align*} Let $\eta \in P^k$ satisfies the boundary condition $D_{+}$. That is, 
\[
0=(1-\frac{1}{n-k+1}L\Lambda)(w_1 \wedge \eta),
\] which is equivalent to $\beta_2|_{\partial M}=0$ and $\beta_4|_{\partial M}=0$. Let $\pi_0: \Omega^k \rightarrow P^k$ be the projection. Then $\partial_{+}=\pi_0 \circ d$ and we get 
\[
\partial_{+}(\rho \partial_{+}\eta)|_{\partial M}=\pi_0(d\rho \wedge \partial_{+}\eta)|_{\partial M}=\pi_0(d\rho \wedge d \eta)|_{\partial M}=\pi_0d(d\rho \wedge \eta)=\partial_{+}(d\rho \wedge \eta)|_{\partial M}.
\] Since $d\rho \wedge \eta=w_{12}\wedge \beta_2+w_1\wedge \beta_4$, there is 
\begin{align*}
\partial_{+}(d\rho \wedge \eta)|_{\partial M}&=\partial_{+}(w_{12}\wedge \beta_2-\frac{1}{n-k+1}\underset{i>1}{\sum}{w_{2i-1, 2i}} \wedge w_1\wedge\beta_3+w_1\wedge \beta_4)|_{\partial M}\\
&=(1-LH^{-1}\Lambda)(w_{12}\wedge d\beta_2)|_{\partial M}+(1-LH^{-1}\Lambda)(w_{1}\wedge d\beta_4)|_{\partial M}=0.\\
\end{align*} Here, we use the fact that $w_{1}\wedge d\beta_2|_{\partial M}=0$ and $w_{1}\wedge d\beta_4|_{\partial M}=0$. This is because $\beta_2|_{\partial M}=0$ and $\beta_4|_{\partial M}=0$ imply that their derivatives along $\partial M$ vanish. Therefore, $\eta \in D_{+}$ implies $\partial_{+}\eta \in D_{+}$.

Since $\partial_{-}=-H^{-1} d^{\Lambda}$ on $P^k$ and $D_{-}$ is equivalent to $JN$, we see that $\partial_{-}\eta \in D_{-}$ when $\eta \in D_{-}$. Applying these two properties to $\mathcal{J}\eta$, we then obtain the other two results.
\end{proof}


As is known, the Dirichlet (Neumann) boundary condition is with respect to the outward normal along the boundary. To be precise, let $\eta \in \Omega^{k}$ and $\overrightarrow{\!\!n}$ be the outward normal. Then $\eta \in D$ if an only if $\eta(v_1, \cdots, v_k)=0$ whenever some vector $v_i=c \overrightarrow{\!\!n}$. And $\eta \in N$ if and only if $\eta (v_1, \cdots, v_k)=0$ whenever no $v_i=c \overrightarrow{\!\!n}$. From the definition, we can see that the $JD/JN$ boundary condition is with respect to the $\mathcal{J} \overrightarrow{\!\!n}$ in this sense. Moreover, when the boundary is of contact type, this vector is exactly given by the Reeb vector field.

To see this, let $(M, \omega)$ be a compact symplectic manifold with a smooth boundary $\partial M$ of contact type. Let $\alpha$ be the contact form and $X_{\alpha}$ be the Reeb vector field. Then there exists the symplectization $(\mathbb{R}\times \partial M, d(e^{a}\alpha))$ such that $(-\epsilon, 0] \times \partial M$ can be identified symplectically with a neighborhood of $\partial M$ in $M$. On $(\mathbb{R}\times \partial M, d(e^{a}\alpha))$, there exists an almost complex structure $J$ and a Riemannian metric $g$ such that 
\begin{itemize}
\item $J$ is invariant under the natural action by $\mathbb{R}-$translation,
\item $\partial_a$ is the outward normal, $\mathcal{J}\partial_a=X_{\alpha}$ and $\mathcal{J}X_{\alpha}=-\partial_a$, where $\partial_a$ denotes the unit vector in the $\mathbb{R}-$direction,
\item $J$ preserves the contact structure.
\end{itemize} 
 Therefore, with respect to this compatible triple $(d(e^{a})\alpha, J, g)$, the $J-$Dirichlet (Neumann) boundary conditions are exactly with respect to the Reeb vector field. See, for example \cite{W}, for details about this structure on the contact boundary.

Lastly, in considering the primitive projection, we have another interpration of the boundary conditions $D_{+}$ and $N_{-}$.
 \begin{lemma}\label{4}
 Let $\pi: \Omega^k \rightarrow P^k$ be the primitive projection, for $k \leq n$. Then it reduces 
 \begin{align*}
 &\pi: \Omega^k_D \rightarrow P^k_{D_{+}}, k<n,\\
  &\pi: \Omega^k_{JD} \rightarrow P^k_{N_{-}}, k<n.
 \end{align*}
 \end{lemma}
 \begin{proof} Let $\eta \in \Omega^k$ satisfies the Dirichlet boundary conditions. We write $\eta$ in the following expression:
 \[
 \eta=p^k+\omega \beta
 \] with $p^k \in P^k$ and $\beta \in \Omega^{k-2}$. Thus, $\pi (\eta)=p^k$. And $\eta \in D$ implies that 
\[
 0=d\rho \wedge \eta=d\rho \wedge p^k+\omega \wedge (d\rho \wedge \beta) \, \text{ on $\partial M$}.
\] Then $0=\pi(d\rho \wedge \eta)=\pi (d\rho \wedge p^k)$ on $\partial M$. This means that $p^k \in D_{+}$. The other result is implied by applying this result to $\mathcal{J}\eta$.
  \end{proof}

 \section{Hodge theory for symplectic Laplacians }
In this section, we demonstrate the Hodge theory for the symplectic Laplacians by employing the ellipticity of some boundary value problems (BVPs).  We first recall some results from the elliptic theory.
\subsection{Elliptic boundary value problems}
Given a compact manifold $M$ with smooth boundary $\partial M$, let $E$  be a vector bundle over $M$ and $G_j$ be a vector bundle over $\partial M$, for $j=1, \cdots, J$. Consider the following elliptic BVP:
\begin{equation*}
\left\{\begin{aligned}
&P: C^{\infty}(M,E) \rightarrow C^{\infty}(M,E)\\
&B_j: C^{\infty}(M,E) \rightarrow C^{\infty}(\partial M, G_j), \, j=1, \cdots, J
\end{aligned}\right.
\end{equation*} Here $P$ is of order $2m$ and $B_j$ is of order $m_j$. Then the combined operator $\mathcal{P}=\{P, B_j\}$ is Fredholm:
\begin{equation*}
\mathcal{P}: H^s(M,E)\rightarrow H^{s-2m}(M,F)\oplus H^{s-m_1-\frac{1}{2}}(\partial M, G_1) \oplus \cdots \oplus H^{s-m_J-\frac{1}{2}}(\partial M, G_J).
\end{equation*} Moreover we say $\{P, B_j\}$ is self-adjoint if $P$ is self-adjoint and the following holds: for any $u,v \in C^{\infty}(M,E)$,
\begin{itemize}
\item if $B_j(u)=B_j(v)=0$ for every $j$, then $(Pu, v)=(u, Pv)$;
\item if $B_j(u)=0$ for every $j$, and $(Pu, v)=(u, Pv)$, then $B_j(v)=0$ for every $j$.
\end{itemize}
The next lemma follows from the results of elliptic BVPs. For a general reference, see \cite{H3} and \cite{SM}.\begin{lemma}\label{ell}
For the self-adjoint elliptic BVP $\{P, B_j\}$, the following holds: 
\begin{itemize}
\item The kernel of $\mathcal{P}$, denoted by $\ker \mathcal{P}$,  is finite and smooth.  
\item For any $f \bot \ker \mathcal{P}$ in $H^s(M,E)$, there exists the unique $u \in H^{s+2m}(M,E)$ and $u \bot \ker\mathcal{P}$ such that $Pu=f$ and $B_j(u)=0$ for every $j$.
\item If $f \in H^s(M,E)$ and $Pu=f, B_j(u)=0$ for all $j=1, \cdots, J$, then $u \in H^{s+2m}(M, E)$.
\end{itemize}
\end{lemma}
Based on this lemma, we next show that the weak solutions of self-adjoint elliptic BVPs are actually strong solutions.
\begin{lemma}\label{ws}
Given $f \in L^2(M,E)$.  Let $u \in L^2(M,E)$ satisfy the following:
\[
(u, Pv)=(f, v)
\] for any $v \in C^{\infty}(M,E)$ satisfying $B_j(v)=0$, with $j=1, \cdots, J$. Then $u \in H^{2m}(M,E)$ and 
\[ 
Pu=0, \, B_j(u)=0, j=1, \cdots J.
\]
\end{lemma}
When $f=0$, the lemma implies immediately the following corollary.
\begin{cor}\label{ker}
If $u \in L^2(M,E)$ satisfies $(u, Pv)=0$ for any $v\in C^{\infty}(M,E)$ with $B_j(v)=0$, $j=1, \cdots, J$, then $u \in \ker \mathcal{P}$. That is, $u$ is smooth and $B_j(u)=0$ for $j=1, \cdots, J$.
\end{cor}
The proof of Lemma \ref{ws} is based on the argument given by Schechter in \cite{SM}, where the case for functions is proved.
\begin{proof}[Proof of Lemma \ref{ws}]
Since the space $\ker \mathcal{P}$ is finite-dimensional, there is an $f= f^1+f^2$ with $f^1 \in \ker \mathcal{P}$ and $f^2$ orthogonal to $\ker\mathcal{P}$. By Lemma \ref{ell}, there exists a $w \in H^{2m}(M, E)$ such that $P w=f^2$ and $B_j(w)=0$, for $j=1, \cdots, J$. Then 
\[
(u-w, Pv)=(f^1, v)
\] for any $v \in C^{\infty}(M,E)$ satisfying the boundary condition $B_j(v)=0$, for all $j=1, \cdots, J$. There exists a sequence $w_j \in C^{\infty}(M,E)$ such that $w_j \rightarrow u-w$ in $L^2$ norm, as $j \rightarrow \infty$. Let $w_j=w_j^1+w_j^2$ with $w_j^1 \in \ker \mathcal{P}$ as the projection and $w_j^2 \bot \ker\mathcal{P}$. Then there exists $v_i \in H^{2m}(M,E)$ and $v_i \bot \ker \mathcal{P}$ such that $P v_i=w_i^1$ and $B_j(v_i)=0$ for every $i, j$. Therefore  
\begin{align*}
(u-w, w_j)&=(u-w, w_j^1)+(u-w, w_j^2)=(u-w, Pv_j)+(u-w, w_j^2)\\
&=(f^1, v_j)+(u-w, w_j^2)=(u-w, w_j^2).
\end{align*} As $j\rightarrow \infty$, we get $w_j^2 \rightarrow u-w$. Since $\ker \mathcal{P}$ is closed, $u-w \in \ker \mathcal{P}$ which implies that $ u \in H^{2m}(M,E)$ and $B_j(u)=0$ for all $j=1, \cdots, J$.
\end{proof}

\subsection{Hodge theory for symplectic Laplacians of second-order }

\begin{definition}
Call the following sets of differential forms spaces of harmonic fields for corresponding Laplacians :
 \begin{align*}
P\mathcal{H}_{+}^k&=\{\eta \in H^1P^k|\, \partial_{+} \eta=\partial_{+}^{\ast} \eta=0\},\\
P\mathcal{H}_{-}^k&=\{\eta \in H^1P^k|\,\partial_{-} \eta=\partial_{-}^{\ast} \eta=0\},\\
P\mathcal{H}_{++}^n&=\{\eta \in H^2P^n|\,\partial_{+} \partial_{-}\eta=\partial_{+}^{\ast} \eta=0\},\\
P\mathcal{H}_{--}^n&=\{\eta \in H^2P^n|\, \partial_{-} \eta=\partial_{-}^{\ast} \partial_{+}^{\ast}\eta=0\}.
 \end{align*} 
  \end{definition} 
\begin{remark} The concepts of harmonic fields are different from that of harmonic forms for an operator when the boundary is not vanishing. For example, a form $\eta \in P^k$ is the harmonic form for $\Delta_{+}$ if and only if $\Delta_{+} \eta=0$ on $M$. However, this does not imply that $\eta$ is a harmonic field when the boundary is not vanishing. 
\end{remark} 
The elliptic theory implies the following result for $\Delta_{+}$. (We follow the convention where the additional subscript, e.g. $D_+$ and $N_+$, identifies the boundary condition that the differential forms satisfy.)  
\begin{thm}[Hodge decomposition for $\Delta_{+}$]
For $k<n$, 
\begin{itemize}
\item $P\mathcal{H}^k_{+, D_{+}}$ and $P\mathcal{H}^k_{+, N_{+}}$ are finite-dimensional and smooth;
\item  The following decompositions hold:
\begin{align*} 
L^2P^k&=P\mathcal{H}_{+,D_{+}}^k\oplus \partial_{+}H^1\!P_{D_{+}}^{k-1}\oplus \partial_{+}^{\ast}H^1\!P^{k+1},\\ 
L^2P^k&=P\mathcal{H}_{+,N_{+}}^k\oplus \partial_{+}H^1\!P^{k-1}\oplus \partial_{+}^{\ast}H^1\!P_{N_{+}}^{k+1},\\
L^2 P^k&=L^2P\mathcal{H}_{+}^k \oplus \partial_{+}H^1\!P_{D_{+}}^{k-1}\oplus \partial_{+}^{\ast} H^1 P^{k+1}_{N_{+}}.
\end{align*} 
\end{itemize}
\end{thm} 
Applying the above results to $\mathcal{J}\eta$, we obtain the Hodge decompositions for $\Delta_{-}$.
 \begin{cor}[Hodge decomposition for $\Delta_{-}$]
 For $k<n$,
 \begin{itemize}
  \item $P\mathcal{H}^k_{-, D_{-}}$ and $P\mathcal{H}^k_{-, N_{-}}$ are finite-dimensional and smooth.
\item The following decompositions hold:
\begin{align*}
L^2P^k&=P\mathcal{H}_{-, D_{-}}^k \oplus \partial_{-} H^1\!P_{D_{-}}^{k+1}\oplus  \partial_{-}^{\ast}H^1\!P^{k-1},\\
L^2P^k&=P\mathcal{H}_{-, N_{-}}^k \oplus  \partial_{-}H^1\!P^{k+1}\oplus  \partial_{-}^{\ast}H^1\!P_{N_{-}}^{k-1},\\
L^2 P^k&=L^2P\mathcal{H}_{-}^k \oplus  \partial_{-}H^1 P^{k+1}_{D_{-}}\oplus  \partial_{-}^{\ast}H^1\!P_{N_{-}}^{k-1}.
\end{align*}
\end{itemize}  
\end{cor} To prove this theorem, we first consider the following BVP: for any $\phi, \psi \in P^k$,


\begin{prop} The following boundary value problem is self-adjoint and elliptic for any $\phi, \psi \in P^k$:
\begin{align}\label{BVP1}
&\quad\Delta_{+} \phi=\psi, \, \text{on $M$}\\
&\left\{  \begin{aligned}
& \partial_{+}(\rho\phi)=0, \, \text{on $\partial M$}\\
&\partial_{+}(\rho \partial_{+}^{\ast} \phi)=0, \, \text{on $\partial M$}.
\end{aligned}\right.\nonumber
\end{align}
\end{prop}
\begin{proof}
We first show that this BVP is self-adjoint.  By Green's formula, for any $u, v \in P^k$, there is
\begin{align*}
(\Delta_{+}u, v)&=(\partial_{+}u, \partial_{+}v)+(\partial_{+}^{\ast}u, \partial_{+}^{\ast}v)+\int_{\partial M}\langle \partial_{+}(\rho\partial^{\ast}_{+}u), v\rangle-\langle\partial_{+}u, \partial_{+}(\rho v)\rangle\\
&=(u, \Delta_{+}v)+\int_{\partial M}\langle \partial_{+}(\rho\partial^{\ast}_{+}u), v\rangle-\langle\partial_{+}u, \partial_{+}(\rho v)\rangle+\langle\partial_{+}(\rho u), v\rangle-(u, \partial_{+}(\rho\partial^{\ast}v)).
\end{align*} Thus this BVP is self-adjoint. We know that $\Delta_{+}$ is elliptic on $P^k$. Fix a point $x \in \partial M$ and choose a normal basis $\{w_i\}$ of $\Omega^{\ast}_x$ such that $w_1=d\rho$ and $\omega=w_1\wedge w_2+\cdots+w_{2n-1}\wedge w_{2n}$. Then $\mathcal{J}w_1=w_2$. In order to show that the BVP is elliptic, we need to show that for any $\xi\bot w_1$ in $\Omega^1_x $, if $f(s)$ is a non-increasing solution of 
\begin{equation}\label{s1}
\sigma(\Delta_{+})_x(\xi+iw_1\partial_s)f(s)=0, 
\end{equation} and satisfies 
\[
\sigma(B_i)(\xi+iw_1 \partial_s)f(s)|_{s=0}=0, \, i=1,2, 
\] then $f=0$. Here $B_1(\phi)=\partial_{+}(\rho\phi)$ and $B_2(\phi)=\partial_{+}(\rho \partial_{+}^{\ast}\phi)$ are boundary operators. Let $f(s)=w_1\wedge \beta_1+w_2\wedge \beta_2+(w_1\wedge w_2-\frac{1}{H+1}\underset{i>1}{\sum}w_{2i-1}\wedge w_{2i})\beta_3+\beta_4$. Here, the $\beta_i$'s are primitive forms that are functions of $s$ and have no $w_1$ and $w_2$ components. Without a doubt, it is enough to consider the case of $\xi=w_2$ and $\xi=w_3$. From \ref{s1}, $\beta_i(s)=\beta_i(0)\exp(-c_is)$ for some positive constant $c_i$ in both cases. Moreover, 
\begin{align*}
\sigma(B_1)(\xi+iw_1 \partial_s)f(s)|_{s=0}=0 &\Rightarrow (1-\frac{1}{H}\Lambda)w_1\wedge f(0)=0\\
& \Rightarrow \beta_2(0)=0, \beta_4(0)=0\\
\sigma(B_2)(\xi+iw_1 \partial_s)f(s)|_{s=0}=0 &\Rightarrow (1-\frac{1}{H}\Lambda)w_1\wedge (-i_{w_1}f'(0))=0 \\
&\Rightarrow \beta'_1(0)=0, \beta'_3(0)=0.
\end{align*} Therefore $f(s)=0$ for any $s$. This proves the ellipticity of this BVP.
\end{proof}

We then use Lemma \ref{ell} and Corollary \ref{ker} to obtain the following.
\begin{cor}
For any $k<n$, the space of harmonic fields $P\mathcal{H}^k_{+, D_{+}}$ is finite-dimensional and smooth, and 
\begin{equation}
L^2P^k=P\mathcal{H}^k_{+, D_{+}}\oplus \partial_{+} H^1\!P^{k-1}_{D_{+}}\oplus \partial_{+}^{\ast}H^1\!P^{k+1}.
\end{equation}Moreover, $ \partial_{+} H^1\!P^{k}_{D_{+}}$ is closed in $L^2-$topology for $k<n-1$.
\end{cor}
\begin{proof}
First we are going to show that $P\mathcal{H}^k_{+, D_{+}}$ is the kernel of the BVP \eqref{BVP1}. Let $\eta \in P\mathcal{H}^k_{+, D_{+}}$ and $\phi \in P^k$ with $\phi \in D_{+}$ and $\partial_{+}^{\ast}\phi \in D_{+}$. By Green's formula, there is
\[
0=(\partial_{+}\eta, \partial_{+}\phi)+(\partial_{+}^{\ast}\eta, \partial_{+}^{\ast}\phi)=(\eta, \Delta_{+}\phi)
\]This implies that $\eta$ belongs to the kernel of BVP \eqref{BVP1} by Corollary \ref{ker}. Obviously, the kernel of BVP \eqref{BVP1} is a subset of $P\mathcal{H}^k_{+, D_{+}}$. Therefore, as the kernel of BVP \eqref{BVP1}, $P\mathcal{H}^k_{+, D_{+}}$ is finite-dimensional and smooth by Lemma \ref{ell}.

Thus,
\[
L^2P^k=P\mathcal{H}^k_{+, D_{+}}\oplus P\mathcal{H}_{+, D_{+}}^{k,\bot}
\] where $P\mathcal{H}_{+, D_{+}}^{k,\bot}$ denotes the orthogonal complement. For any $\eta \in L^2P^k$, let $\eta_1$ be its projection to $P\mathcal{H}^k_{+, D_{+}}$. By Lemma \ref{ws}, there exists a unique $\phi \in H^2\!P^k \cap P\mathcal{H}_{+, D_{+}}^{k,\bot}$ that solves \eqref{BVP1} with $\psi=\eta - \eta_1$.
Therefore, we can write 
\[
\eta=\eta_1+\partial_{+}(\partial_{+}^{\ast}\phi)+\partial_{+}^{\ast}(\partial_{+}\phi)
\] with $\eta_1 \in P\mathcal{H}^k_{+, D_{+}}$ and $\partial_{+}^{\ast}\phi\in H^1\!P^{k-1}_{D_{+}}$. This proves the decomposition. Moreover, we also conclude the $L^2-$closeness of $\partial_{+}H^1\!P^k_{D_{+}}$ from this decomposition by the standard functional analysis argument.
\end{proof}
Consider another self-adjoint elliptic BVP: for any $\phi, \psi \in P^k$,
\begin{align*}
&\quad\Delta_{+} \phi=\psi, \, \text{on $M$}\\
&\left\{  \begin{aligned}
& \partial_{+}^{\ast}(\rho\phi)=0, \, \text{on $\partial M$}\\
&\partial_{+}^{\ast}(\rho \partial_{+} \phi)=0, \, \text{on $\partial M$}. 
\end{aligned}\right.
\end{align*}The ellipticity of this BVP implies the following by similar arguments as above.
\begin{cor}
For $k < n$, the space of harmonic fields $P\mathcal{H}^k_{+, N_{+}}$ is finite-dimensional and smooth, and there is the decomposition
\begin{equation*}
L^2P^k=P\mathcal{H}^k_{+, N_{+}}\oplus \partial_{+} H^1\!P^{k-1}\oplus \partial_{+}^{\ast}H^1\!P^{k+1}_{N_{+}}.
\end{equation*}Moreover $ \partial_{+} H^1\!P^{k}_{N_{+}}$ is closed in $L^2-$topology for $k<n$.
\end{cor}
To complete the proof of the theorem, we only need to prove the following decomposition.
\begin{prop}
The following orthogonal decomposition holds:
\[
L^2P^k=P\mathcal{H}^k_{+}\oplus \partial_{+}H^1\!P^{k-1}_{D_{+}}\oplus \partial_{+}^{\ast}H^1\!P_{N_{+}}^{k+1}.
\]
\end{prop} We shall follow similar arguments for $\Delta_d$ in \cite{S} to prove this proposition.
\begin{proof}
For any $\eta \in L^2P^k$, by the corollaries above, we have the decomposition 
\begin{align*}
\eta&=\eta_1+\partial_{+}\alpha_1+\partial_{+}^{\ast}\beta_1\\
\eta&=\eta_2+\partial_{+}\alpha_2+\partial_{+}^{\ast}\beta_2
\end{align*}with $\eta_1 \in P\mathcal{H}^k_{+, D_{+}}, \eta_2\in P\mathcal{H}^{k}_{+, N_{+}}, \alpha_1 \in D_{+}$ and $\beta_2 \in N_{+}$. Let $u=\eta-\partial_{+}\alpha_1-\partial_{+}^{\ast}\beta_2$. We first show that $u \in P\mathcal{H}^k$ when $\eta \in H^1\!P^k$. This is because
\begin{align*}
(u, \partial_{+}v)&=(\eta-\eta_1, \partial_{+}v)-(\partial_{+}\alpha_1, \partial_{+}v)=0\,  \text{for}\,v\in H^1\!P^{k-1}_{D_{+}}\\
(u, \partial_{+}^{\ast}v)&=(\eta-\eta_2, \partial_{+}^{\ast}v)-(\partial_{+}^{\ast}\beta_2, \partial_{+}^{\ast}v)=0\, \text{for}\, v\in H^1\!P^{k+1}_{N_{+}}, 
\end{align*} and $H^1\!P^{k}_{D_{+}}$ and $H^1\!P^{k}_{N_{+}}$ are dense in $H^1\!P^k$. Therefore, we obtain 
\[
H^1\!P^k=P\mathcal{H}^k_{+}\oplus \partial_{+}H^1\!P^{k-1}_{D_{+}}\oplus \partial_{+}^{\ast}H^1\!P_{N_{+}}^{k+1}.
\] Since $ \partial_{+}H^1\!P^{k-1}_{D_{+}}$ and $\partial_{+}^{\ast}H^1\!P_{N_{+}}^{k+1}$ are closed in the $L^2-$topology, the $L^2-$decomposition is clear by means of a completion argument.
\end{proof}

\subsection{Hodge theory for symplectic Laplacians of fourth-order}
Let us define some boundary conditions first. 
\begin{definition}
We say  $\eta \in D_{+-}$ if $\eta \in D_{-}$ and $\,\partial_{-}\eta \in D_{+}$. Further, we say $\eta \in N_{+-}$ if $\eta \in N_{+}$ and $\,\partial_{+}^{\ast}\eta \in N_{-}$. 
\end{definition}
\begin{thm}[Hodge decompositions for $\Delta_{++}$]
Consider $\Delta_{++}$ on $P^n$.  Then,
\begin{itemize}
\item $P\mathcal{H}^n_{++, N_{+}}$ and $P\mathcal{H}^n_{++, D_{+-}}$ are finite-dimensional and smooth; 
\item
The following decompositions hold:
\begin{align*}
L^2P^n&=P\mathcal{H}^n_{++, N_{+}}\oplus \partial_{-}^{\ast}\partial_{+}^{\ast}H^2\!P^n_{N_{+-}}\oplus \partial_{+}H^1\!P^n,\\ 
L^2P^n&=P\mathcal{H}^n_{++, D_{+-}}\oplus  \partial_{-}^{\ast}\partial_{+}^{\ast}H^2\!P^n\oplus \partial_{+}H^1\!P^n_{D_{+}},\\
L^2P^n&=L^2P\mathcal{H}_{++}^n \oplus  \partial_{-}^{\ast}\partial_{+}^{\ast}H^2\!P^n_{N_{+-}}\oplus H^1\!P^n_{D_{+}}.
\end{align*}
\end{itemize}
\end{thm}
We can obtain similar results for $\Delta_{--}$ by applying the above theorem to $\mathcal{J}\eta$.
\begin{cor}[Hodge decompositions for $\Delta_{--}$]
Consider $\Delta_{--}$ on $P^n$.  Then,
\begin{itemize}
\item $P\mathcal{H}^n_{--, D_{-}}$ and $P\mathcal{H}^n_{--, N_{+-}}$ are finite-dimensional and smooth; 
\item
 The following decompositions hold:
\begin{align*}
L^2P^n&=P\mathcal{H}^n_{--, D_{-}}\oplus \partial_{+}\partial_{-}H^2\!P_{D_{+-}}^n\oplus \partial_{-}^{\ast}H^1\!P^{n-1},\\ 
L^2P^n&=P\mathcal{H}^n_{--, N_{+-}}\oplus  \partial_{+}\partial_{-}H^2\!P^n\oplus \partial_{-}^{\ast}H^1\!P^{n-1}_{N_{-}},\\
L^2P^n&=L^2P\mathcal{H}_{--}^n \oplus \partial_{+}\partial_{-}H^2\!P_{D_{+-}}^n\oplus H^1\!P^{n-1}_{N_{-}}.
\end{align*}
\end{itemize}
\end{cor}
Following the same arguments as for the case of $\Delta_{+}$, we prove the theorem for $\Delta_{++}$ using the ellipticity of some BVPs.
\begin{prop}
The following boundary value problem is self-adjoint and elliptic for any $\phi, \psi \in P^n$:
\begin{align}\label{BVP2}
&\quad\Delta_{++}\phi=\psi\, \text{on $M$}\\
&\left\{\begin{aligned}
&\partial_{+}^{\ast}(\rho \phi)=0 \, \text{on $\partial M$}\\
&\partial_{+}^{\ast}(\rho \partial_{+}\partial_{-} \phi)=\partial_{+}^{\ast}(\rho \partial_{+}\partial_{+}^{\ast}\phi)=0\, \text{on $\partial M$}\\
&\partial_{-}^{\ast}(\rho \partial_{+}^{\ast}\partial_{+}\partial_{-}\phi)=0\, \text{on $\partial M$}.
\end{aligned}\right.\nonumber
\end{align}
\end{prop}
\begin{proof}
By Green's formula, we have 
\begin{align*}
(\Delta_{++}u, v)&=(\partial_{+}\partial_{-}u, \partial_{+}\partial_{-}v)+(\partial_{+}\partial_{+}^{\ast}u,\partial_{+}\partial_{+}^{\ast}v )\\
&+\int_{\partial M}\Big\{\langle\partial_{-}^{\ast}(\rho \partial_{+}^{\ast}\partial_{+}\partial_{-}u), v\rangle+\langle\partial_{+}^{\ast}(\rho\partial_{+}\partial_{-}u), \partial_{-}v\rangle\\
&\qquad\quad -\langle \partial_{+}^{\ast}\partial_{+}\partial_{+}^{\ast}u, \partial_{+}^{\ast}(\rho v)\rangle+\langle \partial_{+}^{\ast}(\rho \partial_{+}\partial_{+}^{\ast}u), \partial_{+}^{\ast}v\rangle\Big\}.
\end{align*} Therefore, 
\[
(\Delta_{++}u, v)=(\partial_{+}\partial_{-}u, \partial_{+}\partial_{-}v)+(\partial_{+}\partial_{+}^{\ast}u,\partial_{+}\partial_{+}^{\ast}u)
\] whenever $u$ and  $v$ satisfy the boundary conditions
\begin{align*}
\partial_{+}^{\ast}(\rho v)|_{\partial M}&=0\\
\partial_{+}^{\ast}(\rho \partial_{+}\partial_{-} u)|_{\partial M}&=\partial_{+}^{\ast}(\rho \partial_{+}\partial_{+}^{\ast}u)|_{\partial M}= \partial_{-}^{\ast}(\rho \partial_{+}^{\ast}\partial_{+}\partial_{-}u)|_{\partial M}=0.
\end{align*} This fact implies that the BVP is self-adjoint.  Like the proof above, we choose the basis $\{w_i\}$. Write $f(s)=w_1\wedge \beta_1+w_2\wedge\beta_2+(w_1\wedge w_2-\underset{i>1}{\sum}w_{2i-1}w_{2i})\beta_3$. Here again, $w_i$'s are primitve forms containing neither $w_1$ nor $w_2$ and are functions of $s$. The equation 
\[
\sigma(\Delta_{++})(\xi+iw_1\partial_s)f(s)=0
\] then implies that $\beta_i(s)=\beta_i(0)\exp(-c_is)$ for some positive constant $c_i$. Let $B_1(\eta)=\partial_{+}^{\ast}(\rho \phi), B_2(\eta)=\partial^{\ast}_{+}(\rho \partial_{+}\partial_{-}\eta)$, then 
\begin{align*}
\sigma(B_1)(\xi+iw_1\partial_s)\eta(0)=0&\Rightarrow \beta_1(0)=0, \beta_3(0)=0\\
\sigma(B_2)(\xi+iw_1\partial_s)\eta(0)=0&\Rightarrow i_{w_1}\left(\sigma(\partial_{+}\partial_{-})(\xi)\eta-\sigma(\partial_{+}\partial_{-})(w_1)\eta''\right)|_{s=0}=0\Rightarrow \beta_2'(0)=0.
\end{align*}Therefore $f(s)=0$ for any $s$, and thus, the ellipticity of this BVP follows.
\end{proof}
The ellipticity of \eqref{BVP2} yields the following.
\begin{cor}
The space $\mathcal{H}^n_{++, N_{+}}$ is finite-dimensional and smooth. And there is the decomposition
\[
L^2P^n=P\mathcal{H}^n_{++, N_{+}}\oplus \partial_{-}^{\ast}\partial_{+}^{\ast}H^2\!P_{N_{+-}}^n\oplus \partial_{+}H^1\!P^n.\] 
Moreover, the space $\partial_{-}^{\ast}\partial_{+}^{\ast}H^2\!P_{N_{+-}}^n$ is closed in the $L^2-$topology.
\end{cor}
\begin{proof}
By Green's formula, it is not hard to see that 
\[
0=(\partial_{+}\partial_{-}u, \partial_{+}\partial_{-}v)+(\partial_{+}\partial_{+}^{\ast}u, \partial_{+}\partial_{+}^{\ast}v)=(u, \Delta_{++}v)
\] where $u \in P\mathcal{H}^n_{++, N_{+}} $and $v \in P^n$ satisfies the boundary conditions in BVP \eqref{BVP2}. Then Corollary \ref{ws} implies that $u$ belongs to the kernel of the BVP \eqref{BVP2}. Since the kernel of BVP \eqref{BVP2} is a subset of $P\mathcal{H}^n_{++, N_{+}}$ obviously, $P\mathcal{H}^n_{++, N_{+}}$ is the kernel,  and then finite-dimensional and smooth, by Lemma \ref{ell}.

For any $\eta \in L^2P^n$, let $\eta_1$ be its projection to $P\mathcal{H}^n_{++, N_{+}}$. By Lemma \ref{ell}, there exists a unique $\phi \in H^4 P^n$ that solves \eqref{BVP2} with $\psi = \eta - \eta_1$.
Therefore, we obtain
\[
\eta=\eta_1+\partial_{-}^{\ast}\partial_{+}^{\ast}(\partial_{+}\partial_{-}\phi)+\partial_{+}(\partial_{+}^{\ast}\partial_{+}\partial_{+}^{\ast}\phi)
\] with $\eta_1 \in P\mathcal{H}^n_{++, N_{+}}$ and $\partial_{+}\partial_{-}\phi \in H^2P^n_{N_{+-}}$. Moreover, the $L^2$-closedness of $\partial_{-}^{\ast}\partial_{+}^{\ast}H^2\!P_{N_{+-}}^n$ is guaranteed by this decomposition. 
\end{proof} Consider another BVP which is also elliptic and self-adjoint: For $\phi, \psi \in P^n$
\begin{align*}
&\quad\Delta_{++}\phi=\psi\, \text{on $M$}\\
&\left\{\begin{aligned}&\partial_{-}(\rho \phi)=0 \, \text{on $\partial M$}\\
&\partial_{+}(\rho \partial_{-} \phi)=0\, \text{on $\partial M$}\\
&\partial_{+}^{\ast}(\rho \partial_{+}\partial_{+}^{\ast}\phi)=0\, \text{on $\partial M$}\\
&\partial_{+}(\rho \partial_{+}^{\ast}\partial_{+}\partial_{+}^{\ast}\phi)=0\, \text{on $\partial M$}.
\end{aligned}\right.
\end{align*} By similar argument as above, the ellipticity of this BVP implies the following:
\begin{cor}
The space $P\mathcal{H}^n_{++, D_{+-}}$ is finite-dimensional and smooth. And there is the decomposition 
\[
L^2P^n=P\mathcal{H}^n_{++, D_{+-}} \oplus \partial_{-}^{\ast}\partial_{+}^{\ast} H^2\!P^n\oplus \partial_{+}H^2\!P^{n-1}_{D_{+}}. 
\]Moreover, the space $\partial_{+}H^2\!P^{n-1}_{D_{+}}$ is closed in the $L^2-$topology.
\end{cor}
And similar to arguments given for the case of $\Delta_{+}$, two proceeding corollaries above together imply the following decomposition and complete the proof of this theorem.
\[
L^2P^n=P\mathcal{H}^n_{++} \oplus \partial_{-}^{\ast}\partial_{+}^{\ast} H^2\!P_{N_{+-}}^n\oplus \partial_{+}H^2\!P^{n-1}_{D_{+}}. 
\]
\section{Symplectic cohomology}
In this section, we apply the results obtained in the previous section to study certain cohomologies on compact symplectic manifolds with boundary. Through the isomorphisms we build between cohomologies and harmonic fields, we can demonstrate the finiteness of these symplectic cohomologies. Moreover, these isomorphism imply that the dimensions of spaces of harmonic fields with certain boundary conditions are indeed symplectic invariants.
\subsection{Primitive cohomologies }Recall the symplectic elliptic complex of Section 2:
$$\begin{CD}
0@>\partial_{+}>>P^0@>\partial_{+}>>P^1@>\partial_{+}>>\cdots @>\partial_{+}>>P^{n-1}@>\partial_{+}>>P^n\\
    @. @. @. @.   @.                                                                                                            @VV{\partial_{+}\partial_{-}}V\\
0@<\partial_{-}<<P^0@<\partial_{-}<<P^1@<\partial_{-}<<\cdots @<\partial_{-}<<P^{n-1}@<\partial_{-}<<P^n
\end{CD}$$ Tseng and Yau studied the cohomologies of this complex in \cite{TY2}, which we shall write as follows:
\begin{align*}
PH^k(\partial_{+})&=\frac{\ker \partial_{+}~|_{P^k}}{\rm{\rm{im}~} \partial_{+}~|_{P^{k-1}}}\, ~~~~\text{for $k<n$},\\
PH^n(\partial_{+})&=\frac{\ker \partial_{+}\partial_{-}~|_{P^n}}{\rm{\rm{im}~} \partial_{+}~|_{P^{n}}},\\
PH^n(\partial_{-})&=\frac{\ker \partial_{-}~|_{P^{n}}}{\rm{\rm{im}~} \partial_{+}\partial_{-}~|_{P^{n}}},\\
PH^k(\partial_{-})&=\frac{\ker \partial_{-}~|_{P^{k}}}{\rm{\rm{im}~} \partial_{-}~|_{P^{k+1}}}\, ~~~~\text{for $k<n$}.
\end{align*} 
Through the Hodge decompositions, we find the following properties of these cohomologies on manifolds with boundary. 
\begin{thm}
Let $(M, \omega)$ be a compact symplectic manifold with a smooth boundary. Let $(\omega, J, g)$ be a compatible triple on $M$. Then there are isomorphisms: 
 \begin{align*}
PH^k(\partial_{+})&\cong P\mathcal{H}^k_{+, N_{+}},\, PH^k(\partial_{-})\cong P\mathcal{H}^k_{-,N_{-}} \\
PH^n(\partial_{+})&\cong P\mathcal{H}^n_{++, N_{+}},\, PH^n(\partial_{-})\cong P\mathcal{H}^n_{--,N_{+-}}\\
\end{align*}
\end{thm}
\begin{proof}Consider the decomposition
\begin{equation*}
P^k=P\mathcal{H}_{+, N_{+}}^k\oplus \partial_{+}P^{k-1}\oplus \partial_{+}^{\ast}P_N^{k+1}.
\end{equation*} For any $\eta \in \ker \partial_{+}|_{P^k}$, we have : 
\[
\eta=  k_{N_{+}} + \partial_{+} \alpha +\partial_{+}^{\ast} \beta 
\] with $k_{N_{+}} \in P\mathcal{H}^{k}_{+, N_{+}}$, $\alpha \in P^{k-1}$ and $\beta \in P_N^{k+1}$. Then the map
\[
PH^k(\partial_{+}) \rightarrow P\mathcal{H}^{k}_{+, JD}: [\eta] \rightarrow k_{N_{+}}
\] is well-defined and isomorphic. This is because that for any $\eta \in \ker \partial_{+}$, the decomposition above yields an unique expression:
\[
\eta= \partial_{+} \alpha+ k_{N_{+}}.
\] This is an easy conclusion of the Green's formula. Similarly, the decompositions
\begin{align*}
P^k&=P\mathcal{H}_{-,N_{-}}^k\oplus \partial_{-}P^{k+1}\oplus \partial_{-}^{\ast}P_{N_{-}}^{k-1}\\
P^n&=P\mathcal{H}_{++, N_{+}}^n\oplus \partial_{+}P^{n-1}\oplus \partial_{-}^{\ast}\partial_{+}^{\ast}P_{N_{+-}}^{n}\\
P^n&=P\mathcal{H}_{--, N_{+-}}^n\oplus \partial_{+}\partial_{-}P^{n}\oplus \partial_{-}^{\ast}P^{n-1}_{N_{-}}
\end{align*} imply the other isomorphisms of this theorem.
\end{proof} This theorem does not only tell us the finiteness of these primitive cohomologies, but also imply that the dimensions of the spaces of harmonic fields appearing here can be regarded as symplectic invariants. For the spaces of harmonic fields with other boundary conditions, we also find isomorphisms between them and cohomologies. We consider the dual complex: 
$$\begin{CD}
0@<\partial_{+}^{\ast}<<P^0@<\partial_{+}^{\ast}<<P^1@<\partial_{+}^{\ast}<<\cdots @<\partial_{+}^{\ast}<<P^{n-1}@<\partial_{+}^{\ast}<<P^n\\
    @. @. @. @.   @.                                                                                                            @AA{\partial_{-}^{\ast}\partial_{+}^{\ast}}A\\
0@>\partial_{-}^{\ast}>>P^0@>\partial_{-}^{\ast}>>P^1@>\partial_{-}^{\ast}>>\cdots @>\partial_{-}^{\ast}>>P^{n-1}@>\partial_{-}^{\ast}>>P^n.
\end{CD}$$ Denote the corresponding cohomologies of this dual complex by $PH^k(\partial^{\ast}_{+})$ for the upper level and $PH^k(\partial^{\ast}_{-})$ for the lower level, with $k \leq n$. Then the following isomorphisms hold.
\begin{thm}
With the assumption as above, we have the following isomorphisms:
 \begin{align*} 
 PH^k(\partial^{\ast}_{+})& \cong P\mathcal{H}^k_{+,D_{+}}, \, PH^k(\partial_{-}^{\ast})\cong P\mathcal{H}^k_{-, D_{-}}\\
 PH^n(\partial^{\ast}_{+})&\cong P\mathcal{H}^n_{++, D_{+-}},\, PH^n(\partial^{\ast}_{-})\cong P\mathcal{H}^n_{--, D_{-}}.
\end{align*}
\end{thm}
For spaces of harmonic fields appearing in this theorem, their dimensions are also symplectic invariants. This is because the following isomorphisms induced by the operator $\mathcal{J}$:
\begin{align*}
P\mathcal{H}^{k}_{+, D_{+}}&\cong P\mathcal{H}^{k}_{-, N_{-}}, \, P\mathcal{H}^{k}_{+, D_{-}}\cong P\mathcal{H}^{k}_{-, N_{+}}\\
P\mathcal{H}^n_{++, D_{+-}} & \cong P\mathcal{H}^{n}_{-+, D_{+-}}, \, P\mathcal{H}^n_{++, N_{+}} \cong P\mathcal{H}^{n}_{-+, D_{-}}.
\end{align*} Moreover, the operator $\mathcal{J}$ reduce the isomorphisms between these symplectic cohomologies and the dual cohomologies:
\[
 PH^k(\partial_{+})\cong PH^k(\partial^{\ast}_{-}),\, PH^k{(\partial_{+}^{\ast}})\cong PH^k(\partial_{-}).
\] for $k \leq n$.

\subsection{Relative symplectic cohomologies}
Like the de Rham case, we can talk about relative symplectic cohomologies by posing suitable boundary conditions on the primitive complex above. In fact,  we obtain:
$$\begin{CD}\label{rc}
0@>\partial_{+}>>P_{D_{+}}^0@>\partial_{+}>>P_{D_{+}}^1@>\partial_{+}>>\cdots @>\partial_{+}>>P_{D_{+}}^{n-1}@>\partial_{+}>>P_B^n\\
    @. @. @. @.   @.                                                                                                            @VV{\partial_{+}\partial_{-}}V\\
0@<\partial_{-}<<P^0@<\partial_{-}<<P^1_{D_{-}}@<\partial_{-}<<\cdots @<\partial_{-}<<P_{D_{-}}^{n-1}@<\partial_{-}<<P^n_{D_{-}}
\end{CD}$$ with $P^n_B=\{\eta \in P^n~|~\partial_{-}(\rho\partial_{+}\partial_{-}\eta)|_{\partial m})=0\}$. This complex is well-defined, since $\partial_{+}$ preserved the boundary condition $D_{+}$ and $\partial_{-}$ preserved $D_{-}$. The corresponding cohomology of this complex is denoted by $PH^k(\partial_{+}, D_{+})$ for the upper level and $PH^k(\partial_{-}, D_{-})$ for the lower level, for $k \leq n$. The usual boundary conditions like the Dirichlet or Neumann, do not yield relative cohomologies in this case since they are not preserved by the operators in this complex.

By using the decompositions we obtained, these relative cohomologies are also isomorphic to spaces of harmonic fields with certain boundary conditions.
\begin{thm}
Given the same assumption as above, we have the following isomorphisms for $k<n$, 
\[
PH^k(\partial_{+}, D_{+}) \cong P\mathcal{H}^k_{+, D_{+}}, \, PH^k(\partial_{-}, D_{-}) \cong P\mathcal{H}^k_{+, D_{-}}.
\]
\end{thm}
\begin{proof} The argument is similar to the one above.  Let us point out the decompositions that imply the corresponding isomorphisms. The decomposition
\[
P^k=P\mathcal{H}^k_{D_{+}} \oplus \partial_{+}P^{k-1}_{D_{+}} \oplus \partial_{+}^{\ast}P^{k+1}.
\] implies the isomorphism
\[
PH^k(\partial_{+}, D_{+}) \rightarrow P\mathcal{H}^k_{D_{+}}: [\eta] \rightarrow \lambda.
\]  The decomposition $P^k=P\mathcal{H}^k_{D_{-}} \oplus \partial_{-}P^{k+1}_{D_{-}} \oplus \partial_{-}^{\ast}P^{k-1}$ yields the other isomorphism.
\end{proof}
\begin{remark}
Notice that the same argument does not give an isomorphism between $PH^n(\partial_{+}, D_{+})$ and $P\mathcal{H}^n_{++, D_{+-}}$. In fact, employing the decomposition
\[
P^n=P\mathcal{H}^n_{++, D_{+-}}\oplus \partial_{-}^{\ast}\partial_{+}^{\ast} P^n\oplus \partial_{+}P_{D_{+}}^{n-1}
\] will give a surjective map from $PH^n(\partial_{+}, D_{+})$ to $P\mathcal{H}^n_{++, D_{+-}}$. The other decompositions do not work in this case. The main obstruction is that the boundary condition $P^n_B$ is much weaker than $D_{+-}$. In fact, for the case of $PH^n(\partial_{+}, D_{-})$, no such well-defined map can be found through the decompositions we obtained. 
\end{remark}

We can also consider the dual of the complex above: 
$$\begin{CD}
0@>\partial_{-}^{\ast}>>P^0@>\partial_{-}^{\ast}>>P_{N_{-}}^1@>\partial_{-}^{\ast}>>\cdots @>\partial_{-}^{\ast}>>p_{N_{-}}^{n-1}@>\partial_{-}^{\ast}>>P_C^n\\
    @. @. @. @.   @.                                                                                                            @VV{\partial_{-}^{\ast}\partial_{+}^{\ast}}V\\
0@<\partial_{+}^{\ast}<<P^0@<\partial_{-}<<P^1_{N_{+}}@<\partial_{+}^{\ast}<<\cdots @<\partial_{+}^{\ast}<<P_{N_{+}}^{n-1}@<\partial_{+}^{\ast}<<P^n_{N_{+}}
\end{CD}$$ with $P^n_C=\{\eta\in p^N~|~ \partial_{-}^{\ast}\partial_{+}^{\ast}\eta \in N_{+} \}$. The corresponding cohomologies are denoted by $PH^k(\partial_{+}^{\ast}, N_{+})$ and  $PH^k(\partial_{-}^{\ast}, N_{-})$ for $k \leq n$. Then the operator $\mathcal{J}$ induces the following isomorphisms:
\begin{align*}
&\mathcal{J}: PH^k(\partial_{+}^{\ast}, N_{+})\rightarrow PH^{k}(\partial_{-}, D_{-}), [\eta] \rightarrow [\mathcal{J}\eta],\\
&\mathcal{J}: PH^k(\partial_{-}^{\ast}, N_{-})\rightarrow PH^{k}(\partial_{+}, D_{+}), [\eta] \rightarrow [\mathcal{J}\eta].
\end{align*} 
Combining the isomorphisms above, we obtain the following characterization of these cohomologies
\[
PH^k(\partial_{+}^{\ast}, N_{+})  \cong P\mathcal{H}^k_{+, N_{+}}, \, PH^k(\partial_{-}, N_{-}) \cong P\mathcal{H}^k_{+, N_{-}}\]
for $k<n$.

\section{Boundary value problems}
Another application of the Hodge theory in the boundary case is to solve boundary value problems.  We shall begin with Poincar\'{e} lemmas.
\subsection{Poincar\'{e} lemmas}
\begin{lemma}[Poincar\'{e} lemma for $\partial_{+}$]
Given $(M^{2n}, \omega)$ as a compact symplectic manifold with smooth boundary, let $(\omega, J, g)$ be a compatible triple on it.  For $\eta \in P^k$,  $\eta$ is $\partial_{+}$-exact, i.e. there exists a solution $\phi \in P^{k-1}$ of the equation 
\begin{equation*}
\partial_{+}\phi=\eta
\end{equation*}  if and only if $\eta$ obeys the integrability conditions:
\begin{itemize}
\item when $k<n$
\begin{equation*}
\partial_{+}\eta=0\, \text{and}\, (\eta, \lambda)=0\, \text{ for any $\lambda \in P\mathcal{H}_{+, N_{+}}^k$}.
\end{equation*} 
\item when $k=n$
\begin{equation*}
\partial_{+}\partial_{-}\eta=0\, \text{and}\, (\eta, \lambda)=0\, \text{ $\lambda \in P\mathcal{H}_{++, N_{+}}^n$}.
\end{equation*} 
\end{itemize}
\end{lemma}
\begin{proof}
Consider the Hodge decompositions:
\begin{align*}
P^k&=P\mathcal{H}_{+, N_{+}}^k \oplus \partial_{+}P^{k-1}\oplus \partial_{+}^{\ast}P^{k+1}_{N_{+}}, \, \text{ $k<n$}\\
P^n&=P\mathcal{H}_{++, N_{+}}^n\oplus \partial_{-}^{\ast}\partial_{+}^{\ast}P^n_{N_{+-}}\oplus \partial_{+}P^{n-1}.
\end{align*}  For any $\eta \in P^k$ with $k<n$, it is obvious that $\eta$ satisfies the integrability conditions when $\eta= \partial_{+} \phi$.  Conversely, assume $\eta$ satisfies the integrability conditions. By the decompositions above, there exist $\lambda \in P\mathcal{H}^k_{+, N_{+}}, \alpha \in P^{k-1}$ and $\beta \in P_{N_{+}}^{k+1}$ such that
\[
\eta=\lambda+\partial_{+}\alpha+\partial_{+}^{\ast}\beta.
\] Then 
\begin{equation*}
0=(\partial_{+}\partial_{+}^{\ast}\beta, \beta)=(\partial_{+}^{\ast}\beta, \partial_{+}^{\ast}\beta)
\end{equation*}which implies that $\partial_{+}^{\ast}\beta=0$. Moreover, the integrability condition $(\eta, \lambda)=0$ implies that $\lambda=0$ since $(\eta, \lambda)=(\lambda, \lambda)$. Therefore $\eta=\partial_{+}\alpha$. Similar argument works for the case $k=n$ using the decomposition of $P^n$ above.

\end{proof}
Similarly, we obtain the Poincar\'{e} lemmas for the other differential operators.
\begin{lemma}[Poincar\'{e} lemma for $\partial^{\ast}_{+}$]
 For $k<n$,  $\eta \in P^k$ is $\partial^{\ast}_{+}$ exact, i.e. there exists a solution $\phi \in P^{k+1}$ of the equation 
\begin{equation*}
\partial_{+}^{\ast}\phi=\eta
\end{equation*}  if and only if $\eta$ obeys the integrability conditions:
\begin{equation*}
\partial_{+}^{\ast}\eta=0\, \text{and}\,~ (\eta, \lambda)=0\, \text{ for any $\lambda \in P\mathcal{H}_{+, D_{+}}^k$}.
\end{equation*} 

\end{lemma}

\begin{lemma}[Poincar\'{e} lemma for $\partial_{-}$]
 For $k<n$, $\eta \in P^k$ is $\partial_{-}$ exact, i.e. there exists a solution $\phi\in P^{k+1}$ of the equation 
\begin{equation*}
\partial_{-}\phi=\eta
\end{equation*}  if and only if $\eta$ obeys the integrability conditions:
\begin{equation*}
\partial_{-}\eta=0\, \text{and}\, (\eta, \lambda)=0\, \text{ for any $\lambda \in P\mathcal{H}_{-, N_{-}}^k$}.
\end{equation*} 
\end{lemma}

\begin{lemma}[Poincar\'{e} lemma for $\partial_{-}^{\ast}$]
 A form $\eta \in P^k$ is  is $\partial_{-}^{\ast}$ exact, i.e. there exists a solution $\phi\in P^{k-1}$ of the equation 
\begin{equation*}
\partial_{-}^{\ast}\phi=\eta
\end{equation*}  if and only if $\eta$ obeys the integrability conditions:
\begin{itemize}
\item when $k<n$
\begin{equation*}
\partial_{-}^{\ast}\eta=0\, \text{and}\, (\eta, \lambda)=0\, \text{ for any $\lambda \in P\mathcal{H}_{-, D_{-}}^k$}.
\end{equation*} 
\item when $k=n$
\begin{equation*}
\partial_{-}^{\ast}\partial_{+}^{\ast}\eta=0\, \text{and}\, (\eta, \lambda)=0\, \text{ $\lambda \in P\mathcal{H}_{--, D_{-}}^n$}
\end{equation*} 
\end{itemize}
\end{lemma}

For second order differential operators, we get following results.

\begin{lemma}[Poincar\'{e} lemma for $\partial^{\ast}_{-}\partial_{+}^{\ast}$ and $\partial_{+}\partial_{-}$]
 A form $\eta \in P^n$ is $\partial^{\ast}_{-}\partial_{+}^{\ast}$ exact, i.e. there exists a solution $\phi \in  P^{n}$ of the equation 
\begin{equation*}
\partial_{-}^{\ast}\partial_{+}^{\ast}\phi=\eta
\end{equation*}  if and only if $\eta$ obeys the integrability conditions:
\begin{equation*}
\partial_{+}^{\ast}\eta=0\, \text{and}\, (\eta, \lambda)=0\, \text{ for any $\lambda \in P\mathcal{H}_{++, N_{+-}}^n$}.
\end{equation*} 
A form $\eta \in P^n$ is $\partial_{+}\partial_{-}$ exact, i.e. there exists a solution $\phi \in  P^{n}$ of the equation 
\begin{equation*}
\partial_{+}\partial_{-}\phi=\eta
\end{equation*}  if and only if $\eta$ obeys the integrability conditions:
\begin{equation*}
\partial_{-}\eta=0\, \text{and}\, (\eta, \lambda)=0\, \text{ for any $\lambda \in P\mathcal{H}_{--, D_{+-}}^n$}.
\end{equation*} 
\end{lemma}

\subsection{Harmonic fields and Boundary Value Problems} Besides these Poincar\'{e} lemmas,  there are various BVPs in the symplectic case that can be solved by applying the Hodge decompositions. Here, we consider some BVPs that are related to the existence of harmonic fields. In fact, we will use these BVPs to show that the spaces of harmonic fields are infinite-dimensional without boundary conditions.

\begin{thm}\label{bvpi}
For $\eta \in P^k$  there exists a solution $\phi \in P^{k-1}$of the boundary value problem 
\begin{equation*}
\partial_{+}\phi=\eta\,\, \text{on $M$    and     }\,\partial_{+}(\rho \phi)=\partial_{+}(\rho x)\,\, \text{on $\partial M$}
\end{equation*} with $x \in P^{k-1}$  if and only if $\eta$ and $x$ obey the integrability conditions:
\begin{itemize}
\item when $k<n$
\begin{equation*}
\partial_{+} \eta=0\, \text{and}\, ~ (\eta, \lambda)=\int_{\partial M}\langle \partial_{+}(\rho x), \lambda\rangle
\end{equation*} for any $\lambda \in P\mathcal{H}_{+}^k$.
\item when $k=n$
\begin{equation*}
\partial_{+}\partial_{-}\eta=0\, \text{and}\,~ (\eta, \lambda)=\int_{\partial M}\langle \partial_{+}(\rho x), \lambda \rangle
\end{equation*} for any $\lambda \in P\mathcal{H}^n_{++}$.
\end{itemize}Moreover, the solution $\phi$ can be chosen to satisfy $\partial_{+}^{\ast} \phi=0$.
\end{thm}
\begin{proof}
Given $\eta \in P^k$ and $x \in P^{k-1}$, if there exists $\phi \in P^{k-1}$ such that 
\begin{equation*}
\partial_{+}\phi=\eta\,\, \text{on $M$    and     }\,\partial_{+}(\rho \phi)=\partial_{+}(\rho x)\,\, \text{on $\partial M$},
\end{equation*} it is obvious to see that $\eta$ satisfies the integrability conditions. Now assume $\eta$ satisfies the integrability conditions. When $k<n$, there is the decomposition 
\[
P^k=P\mathcal{H}^k_{+}\oplus \partial_{+} P_{D_{+}}^{k-1}\oplus \partial_{+}^{\ast} P^{k+1}_{N_{+}}. 
\] That is, there exist smooth forms $\lambda \in P\mathcal{H}^k_{+}, \alpha \in P^{k-1}_{D_{+}}$ and $\beta \in P^{k+1}_{N_{+}}$ such that 
\begin{equation*}
\eta=\lambda+\partial_{+}\alpha+\partial_{+}^{\ast}\beta.
\end{equation*} The integrability condition $\partial_{+}\eta=0$ implies that $\partial_{+}\partial_{+}^{\ast}\beta=0$. It then follows that
\[0=(\partial_{+}\partial_{+}^{\ast}\beta, \beta)=(\partial_{+}^{\ast}\beta, \partial_{+}^{\ast}\beta).\]Then $\partial_{+}^{\ast}\beta=0$ and $\eta=\lambda+\partial_{+}\alpha$.
In fact, we can choose $\psi$ to be primitive such that $\partial_{+}(\rho \psi)=\partial_{+}(\rho x)$ on $\partial M$ and $\partial_{+}^{\ast}\psi=0$. By the decomposition above,  let $\tilde{\eta}=\partial_{+}\psi$ and we get
\[
\tilde{\eta}=\tilde{\lambda}+\partial_{+}\tilde{\alpha}
\] with $\tilde{\lambda} \in P\mathcal{H}_{+}^k$ and $\tilde{\alpha} \in P^{k-1}_{D_{+}}$. Now let $\phi=\alpha+\psi-\tilde{\alpha}$ and there is
\begin{align*}
&\partial_{+}\phi=\eta+\tilde{\lambda}-\lambda,\\
&\partial_{+}(\rho \phi)~|_{\partial{M}}=x.
\end{align*} Since $(\tilde{\lambda}-\lambda, \mu)=(\partial_{+}\phi-\eta, \mu)=0$ for any $\mu \in P\mathcal{H}^k_{+}$, it follows that $\tilde{\eta}_h-\eta_h \in P\mathcal{H}^{k, \bot}_{+}$. Thus $\tilde{\eta}_h-\eta_h=0$ and $\phi$ is the solution for this boundary value problem. Via Hodge decompositions,  $\alpha$ and $\tilde{\alpha}$ can be chosen to be $\partial_{+}^{\ast}$-closed. Therefore, so is $\phi$. Similar argument works for the case of $k=n$ using the following decomposition:
\[
P^n=P\mathcal{H}^n_{++}\oplus \partial_{-}^{\ast}\partial_{+}^{\ast}P^n_{N_{+-}}\oplus \partial_{+}P^{n-1}_{D_{+}}.
\]
\end{proof}
\begin{thm}\label{inf}
On a compact symplectic manifold $(M, \omega, J, g)$ with smooth boundary, the space $P\mathcal{H}^k_{+}$ is infinite-dimensional for $0<k<n$.
\end{thm}
\begin{proof}
 For $k<n$, define the map
\[
B: P\mathcal{H}^k_{+} \rightarrow \Omega^{k+1}|_{\partial M}: \eta \rightarrow \partial_{+}(\rho \eta)|_{\partial M}.
\]  From the definition of boundary condition $D_{+}$, we see that $\eta \in D_{+}$ if and only if $B(\eta)=0$. Therefore the kernel of this map is exactly $P\mathcal{H}^k_{+, D_{+}}$ which is finite-dimensional.

We claim that the map $B$ is surjective to the space $\partial_{+}(\rho\partial_{+}P^{k-1})|_{\partial M}$. That is, for any $\psi \in  \partial_{+}P^{k-1}$, there is an $\eta \in P\mathcal{H}^k_{+}$ such that
\begin{align*}
&\partial_{+}\eta=0, \partial_{+}^{\ast}\eta=0, \,\text{on $M$}\\
&\partial_{+}(\rho \eta)=\partial_{+}(\rho \psi), \, \text{ on $\partial M$}.
\end{align*} From Theorem \ref{bvpi}, such an $\eta$ exists as long as 
\[
0=\int_{\partial M}\langle\partial_{+}(\rho \psi),\lambda \rangle
\]  for any $\lambda \in P\mathcal{H}^{k+1}_{+}$ when $k+1<n$, or $\lambda \in P\mathcal{H}^n_{++}$. Since $\psi=\partial_{+}u $, there is 
\[
0=(\partial_{+}\psi, \lambda)=\int_{\partial M}\langle\partial_{+}(\rho \psi),\lambda \rangle.
\]  Since the kernel is finite-dimensional, and $\partial_{+}(\rho\partial_{+}P^{k-1})|_{\partial M}$ is infinite-dimensional, $P\mathcal{H}^k_{+}$ is infinite-dimensional.

\end{proof}
We can also consider a similar BVP involving $\partial_{+}^{\ast}$ and get the following result.
 \begin{thm}
For $\eta \in P^k$ with $k<n$,  there exits a solution $\phi \in P^{k+1}$of the boundary value problem  \begin{equation*}
\partial_{+}^{\ast}\phi=\eta\,\, \text{on $M$ and }\, \partial_{+}^{\ast}(\rho \phi)=\partial_{+}^{\ast}(\rho x)\,\, \text{on $\partial M$}
\end{equation*} with $x \in P^{k}$ if and only if 
\begin{equation*}
\partial_{+}^{\ast} \eta=0\, \text{and} \, (\eta, \lambda)=\int_{\partial M}\langle \partial_{+}^{\ast}(\rho x) , \lambda\rangle
\end{equation*}for any $\lambda \in PH_{+}^k$. Moreover, $\phi$ can be chosen to satisfy that 
\[
\partial_{+}\phi=0\, \text{when $k+1<n$};\,  \partial_{+}\partial_{-}\phi=0, \, \text{ when $k+1=n$}.
\]
\end{thm}
Through similar arguments like above, this theorem implies that $P\mathcal{H}_{++}^n$ is also infinite-dimensional. Moreover the conjugate relations imply the infiniteness of $P\mathcal{H}^k_{-}$ and $P\mathcal{H}^n_{--}$ from that of $P\mathcal{H}^k_{+}$ and $P\mathcal{H}^n_{++}$, respectively.

\section{Discussion}


In this paper, we have investigated some natural boundary conditions that arise from the perspective of symplectic Laplacians.  Perhaps especially noteworthy are two, $D_+$ and $D_-$, whose conditions on the boundary are dependent only on the symplectic structure $\om$.  Here, we will briefly discuss in more details the $D_+$ boundary condition and mention some of its relation with the standard Dirchlet boundary condition (D).  As have already been noted, a differential form that satisfies $D$ automatically satisfies $D_+$, as $D_+$ is generally a weaker boundary condition than $D$.  This can be most easily seen in the projection operation that takes a differential form of degree $k\leq n$ to its primitive component:  $\pi: \Omega^k \rightarrow P^k$.  Adding boundary conditions, we have for $k<n$ 
\begin{equation*}\label{proja}\pi: \Omega^k_D \rightarrow P^k_{D_{+}}.
\end{equation*}
From a different perspective, we know that the Dirichlet boundary condition for forms geometrically imposes the vanishing of a form when pull-backed to the boundary.  In the case of $D_+$, using the association of primitive forms with dual co-isotropic spaces as discussed in \cite{TY1}, the geometric intuition should be that the pullback of a primitive $k$-form is not necessariliy zero, but is zero on any co-isotropic subspaces of co-dimension $k$ on the boundary.  

The above property of projection linking $D$ and $D_+$ can further be used to study relations between the relative primitive cohomologies and the relative de Rham cohomologies through the Lefschetz maps as discussed in \cite{TY3}.  We recall that the Lefschetz map of degree $r$ is a map between the de Rham cohomology $H^*(d)$ by the Lefschetz operator:
\[
L^r: {H}^{n-r}(d)\rightarrow  {H}^{n+r}(d), \qquad [\eta] \rightarrow [\omega^k\wedge \eta].
\] 
This map can be restricted to de Rham elements satisfying the Dirichlet boundary conditions:  
\[
L^r: {H}^{n-r}(d,D)\rightarrow  {H}^{n+r}(d,D).
\]
Following similar exact sequence type arguments as in \cite{TY3}, we can obtain the following relations between certain relative primitve cohomologies and relative de Rham cohomologies, i.e. $H^*(d,D)$, via Lefschetz maps.
\begin{prop}
For $k<n$, we have the isomorphism
\begin{align*}
P{H}^k(\dpp,D_{+}) &\cong \ker[L: {H}^{k-1}(d, D) \rightarrow {H}^{k+1}(d, D)]\\
&\quad \oplus \cok[L: {H}^{k-2}(d, D) \rightarrow {H}^{k}(d, D)].
\end{align*}
\end{prop}
One can ask whether this type of relations may extend to the case $k=n$ and for $PH^k(\dpm, D_-)$.  Interestingly, the same line of reasoning seems to break down exactly in the middle fo the symplectic elliptic complex of Section 5.2.  We believe this is somehow related to the lack of an isomorphism for $PH^n(\dpp, D_+)$ as commented upon in Remark 5.4.  In a sense, the difficulty can be pinpointed to the presence of the second-order differential operator $\dpp\dpm$ in the middle of the elliptic complex.

In fact, the boundary conditions we have discussed in this paper have all been local in nature.  The failure of finding an isomorphism for  $PH^n(\dpp, D_+)$ is an indirect indication that the symplectic complex can not be elliptic with local boundary conditions on manifolds with boundary.   This suggests that one should consider global boundary conditions of Atiyah-Patodi-Singer type \cite{APS} for this complex.  Such should lead to a different type of symplectic invariant which we will discuss in a follow up paper.

\appendix

\section{Gaffney inequalities}
In this appendix, we consider the Gaffney inequality for the second-order symplectic Laplacians, $\Delta_+$ and $\Delta_-$.

Recall first the standard Gaffney inequality \cite{M2} for $\Delta_d$ in the Riemannian case, which will imply a Gaffney inequality for $\Delta_{d^{\Lambda}}$ on symplectic manifolds.  Let $(M,g)$ be a compact Riemannian manifold and define the Dirichlet integral for $\Delta$ as
\begin{equation*}
D_d(\eta, \phi)=(d\eta, d\phi)+(d^{\ast}\eta, d^{\ast}\phi).
\end{equation*} 
\begin{thm}[Gaffney inequality for $\Delta$]
Given a compact Riemannian manifold $(M,g)$ with smooth boundary $\partial M$, there exists a constant $C$ depending only on $g$ such that 
\begin{equation*}
D_d((\eta, \eta)+\|\eta\|^2_0\geq C\|\eta\|^2_1
\end{equation*} for any $\eta \in H^1\Omega^k$ satisfying the Dirichlet boundary condition or Neumann boundary condition.
\end{thm} 
Here, $\|\eta\|_0$ and $\|\eta\|_1$ denote the $L^2-$norm and $H^1-$norm, respectively.   For details of the inequality, see \cite{M1}. Now consider a compact symplectic manifold $(M,\omega)$ with a compatible almost complex structure $J$ and a compatible Riemannian metric $g$. Define the Dirichlet integral $D_{d^{\Lambda}}$ for $\Delta_{d^{\Lambda}}$ as:
\begin{equation*}
D_{d^{\Lambda}}(\eta, \phi)=(d^{\Lambda}\eta, d^{\Lambda}\phi)+(d^{\Lambda \ast}\eta, d^{\Lambda \ast}\phi).
\end{equation*} Then it is easy to see that $D_{d^{\Lambda}}(\eta, \phi)=D_d((\mathcal{J}\eta, \mathcal{J}\phi)$ which is implied by the conjugate relations. Then we get the following Gaffney inequality for $D_{d^{\Lambda}}$.
\begin{cor}[Gaffney inequality for $D_{d^{\Lambda}}$]
Given $(M, \omega)$ as a compact symplectic manifold, let $(\omega, J, g)$ be a compatible triple on it. 
 When $M$ has smooth boundary, there exists a constant $C$ depending only on $(\omega, J, g)$ such that 
\begin{equation*}
D_{d^{\Lambda}}(\eta, \eta) +\|\eta\|^2_0 \geq C\|\eta\|^2_1
\end{equation*} for any $\eta \in H^1\Omega^k$ satisfying the boundary condition $JD$ or $JN$.
\end{cor}

Now let $(M, \omega)$ be a compact symplectic manifold with smooth boundary and a compatible triple $(\omega, \mathcal{J}, g)$. 
\begin{definition}
 For $ \Delta_{+}$, and $\Delta_{-}$, we call the following bilinear forms their Dirichlet integrals, respectively:
\begin{align*}
D_{\partial_{+}}(\eta, \eta)&=(\partial_{+}\eta, \partial_{+}\eta)+(\partial_{+}^{\ast} \eta, \partial_{+}^{\ast}\eta),\\
D_{\partial_{-}}(\eta, \eta)&=(\partial_{-}\eta, \partial_{-}\eta)+(\partial_{-}^{\ast} \eta, \partial_{-}^{\ast}\eta).
\end{align*} 
\end{definition}
\begin{remark}
By a Dirichlet integral for a Laplacian, we mean a bilinear form whose kernel is the same as that of the Laplacian, restricted to differential forms with compact support. It is unique when the Laplacian is of second order.
\end{remark}
\begin{thm}[Gaffney's inequalities for $\Delta_{+}$ and $\Delta_{-}$]\label{thm1}
Let $(M, \omega)$ be a compact symplectic manifold with smooth boundary and a compatible triple $(\omega, \mathcal{J}, g)$. Then there exists a constant $C>0$, depending only on $(\omega, \mathcal{J}, g)$, such that
\begin{align*}
 D_{\partial_{+}}(\eta, \eta)+\|\eta\|^2_{0}&\geq C\|\eta\|^2_{1} ,\\
D_{\partial_{-}}(\eta, \eta)+\|\eta\|^2_{0}&\geq C\|\eta\|^2_{1}
\end{align*} for any $ \eta \in P^{k}$ with $k<n$ satisfying the boundary condition $D$ or $JD$. 
\end{thm}
To simplify the calculations, we introduce the Dirichlet integral $D_{\partial_{-}'}$ for $\Delta_{-}'$:
\[
D_{\partial_{-}'}(\eta, \phi)=(\partial_{-}'\eta, \partial_{-}'\phi)+(\partial_{-}^{'\ast}\eta, \partial_{-}^{'\ast}\phi).
\] For any $\eta \in P^k$, it is easy to see that 
$$D_{\partial_{-}}(\eta, \eta)=\frac{1}{(n-k+1)^2}D_{\partial_{-}'}(\eta, \eta)+\frac{2(n-k)+1}{(n-k)^2(n-k+1)^2}(\partial_{-}^{'\ast}\eta, \partial_{-}^{'\ast}\eta).$$ 
Therefore, the Gaffney inequality will hold for $D_{\partial_-}$ if it is true for $D_{\partial_-'}'$. An advantage to using $D_{\partial_-}'$ is that it relates to $D_{\partial_+}$ via the conjugate relation $D_{\partial_{-}'}(\eta, \eta)=D_{\partial_{+}}(\mathcal{J}\eta, \mathcal{J}\eta)$.

In order to prove Theorem \ref{thm1}, we first recall the following lemma from \cite{G}. 
\begin{lemma}\label{key}
Let $(M,g)$ be a compact Riemannian manifold with boundary and $\nabla$ denote the  Levi-Civita connection. Then, for any $\eta \in \Omega^k$, 
\begin{equation*}
(\nabla \eta, \nabla \eta)= D_d((\eta, \eta)+(\mathcal{R}(\eta), \eta)+BT(\eta).
\end{equation*} Here, $\mathcal{R}$ denotes a curvature operator and $BT(\eta)$ is an integral along $\partial M$ which satisfies
\begin{align*}
   BT(\eta) &= (S(\eta), \eta)_{\partial M} \,\text{ if $\eta \in D$ }, \\
      &=(T(\eta), \eta)_{\partial M}\, \text{ if $ \eta \in N$}.
\end{align*} Here $S$ and $T$ are the curvature operators along $\partial M$ and only depend on the second fundamental form. 
\end{lemma}
In regards to the inner product between forms, the following property can be easily shown to hold:
\begin{lemma}\label{IP}
Let $\eta=L^{r}B_s$ and $\phi=L^{p}B_q$ with $B_s$ and $B_q$ as primitive forms. Assume $k=2r+s=2p+q$. Then $\langle \eta, \phi \rangle =0$ if $r \neq p$. Here $\langle, \rangle$ is the inner product induced by metric $g$.
\end{lemma}

It is useful to introduce 
\begin{align*}
D_J(\eta, \eta)&=D_d(\eta, \eta)-\frac{1}{n-k+1}D_{d^{\Lambda}}(\eta, \eta) \\
&=D_d(\eta, \eta)-\frac{1}{n-k+1}D_d(\mathcal{J}\eta, \mathcal{J}\eta).
\end{align*} 
for any $\eta \in \Omega^k$ with $k<n$.  It is related to $D_{\partial_{+}}$ and $D_{\partial_{-}'}$ as follows:
\begin{lemma}
For any $\eta \in P^k$ with $k<n$, there are
\begin{align*}
D_{\partial_{+}}(\eta, \eta)&=D_J(\eta, \eta)+ \frac{1}{n-k+1}(d^{\Lambda\ast}\eta, d^{\Lambda\ast}\eta)\\
D_{\partial_{-}'}(\eta, \eta)&=D_{\partial_{+}}(\mathcal{J} \eta, \mathcal{J} \eta)\\
&=D_{J}(\mathcal{J} \eta, \mathcal{J} \eta)+\frac{1}{n-k+1}(d \eta, d\eta).
\end{align*}
\end{lemma}
\begin{proof}
Since $\partial_{+}=d-LH^{-1}\dpmm$ on $P^k$, Lemma \ref{IP} implies
\[
(\partial_{+}\eta, \partial_{+}\eta)= (d \eta, d \eta)-\frac{1}{n-k+1}(\dpmm\eta, \dpmm\eta),\, \text{ for any $\eta \in P^k$}.
\] Therefore,
\begin{align*}
D_{\partial_{+}}(\eta, \eta)&=(\partial_{+}\eta, \partial_{+}\eta)+(\partial_{+}^{\ast}\eta, \partial_{+}^{\ast}\eta)\\
&=(d \eta, d \eta)-\frac{1}{n-k+1}(\dpmm\eta, \dpmm\eta)+(d^{\ast} \eta, d^{\ast} \eta)\\
&=D_d(\eta, \eta)-\frac{1}{n-k+1}\left((d^{\Lambda}\eta,d^{\Lambda}\eta)+(d^{\Lambda\ast}\eta,d^{\Lambda\ast}\eta)\right)+\frac{1}{n-k+1}(d^{\Lambda\ast}\eta, d^{\Lambda\ast}\eta)\\
&=D_d(\eta, \eta)-\frac{1}{n-k+1}D_{d^{\Lambda}}(\eta, \eta)+\frac{1}{n-k+1}(d^{\Lambda\ast}\eta, d^{\Lambda\ast}\eta),
\end{align*}having used in the third line the relation $\dpmm=-d^{\Lambda}$ on $P^k$.
The conjugate relation $D_{\partial_{-}'}(\eta, \eta)=D_{\partial_{+}}(\mathcal{J}\eta, \mathcal{J}\eta)$ implies the result for $D_{\partial_{-}}$.
\end{proof}
The following is the key lemma for proving Gaffney's inequalities for second-order symplectic Laplacians.
\begin{lemma}\label{cor1}
There exists a constant $C$ depending only on $(\omega, \mathcal{J}, g)$ such that
\begin{equation*}
D_J(\eta, \eta)+\|\eta\|^2_0\geq C\|\eta\|^2_{1}  
\end{equation*} for any $\eta \in \Omega^k$ with $k<n$ satisfying the following boundary condition labelled by $B_1$:
\begin{equation*}B_1:
\left\{\begin{aligned}
\eta &\in D\, \text{or}\, \eta \in N,\\
\eta &\in JD\, \text{or}\,\eta \in JN.
\end{aligned}\right.
\end{equation*} \end{lemma} 
\begin{proof}
Substitute first the following equality from Lemma \ref{key} into $D_{J}$, we find
\begin{equation*}
(\nabla \eta, \nabla \eta)= D_d(\eta, \eta)+(\mathcal{R}(\eta), \eta)+BT(\eta).
\end{equation*}  
Further,
\begin{align*}
&(\nabla \eta, \nabla \eta)-\frac{1}{n-k+1}(\nabla \mathcal{J}\eta, \nabla \mathcal{J}\eta)=D_d(\eta, \eta)-\frac{1}{n-k+1}D_d(\mathcal{J}\eta, \mathcal{J}\eta)\\
&+(\mathcal{R}(\eta), \eta)+BT(\eta)-\frac{1}{n-k+1}(\mathcal{R}(\mathcal{J}\eta), \mathcal{J}\eta)-\frac{1}{n-k+1}BT(\mathcal{J}\eta)
\end{align*} for any $\eta \in P^k$.
We define the term $ET$ such that :
\[
( \nabla\mathcal{J}\eta, \nabla\mathcal{J}\eta)=(\nabla \eta, \nabla \eta)+ET.
\] In fact, $ET$ is an integral only involving $|\eta|$ and $|\eta||\nabla \eta|$ and we obtain
\[
|ET| \leq C(\epsilon)\|\eta\|_{0}+\epsilon \|\nabla \eta\|_{0}
\] with a constant $C(\epsilon)$ depending only on $(\omega, \mathcal{J}, g)$. Next, we derive
\begin{align*}
\frac{n-k}{n-k+1}(\nabla \eta, \nabla \eta)&=D_J(\eta, \eta)+(\mathcal{R}(\eta), \eta)-\frac{1}{n-k+1}(\mathcal{R}(\mathcal{J}\eta), \mathcal{J}\eta)\\
&+\frac{1}{n-k+1}ET+BT(\eta)-\frac{1}{n-k+1}BT(\mathcal{J}\eta).
\end{align*}
We need to estimate the right hand side of the above equation. For the curvature terms, the compactness of $M$ implies that 
\[
|(\mathcal{R}\eta, \eta)| \leq \|\eta\|_{0}\,\|\mathcal{R}\eta\|_{0}\leq C_{\mathcal{R}}\|\eta\|_{0}.
\] Here, $C_{\mathcal{R}}$ is given by the operator norm of $\mathcal{R} \in End(\Lambda^{\ast})$.  Moreover
\[
|(\mathcal{R}\mathcal{J}\eta, \mathcal{J}\eta)| \leq C_{\mathcal{R}}\|\mathcal{J}\eta\|_{0}=C_{\mathcal{R}}\|\eta\|_{0}.\]For the boundary terms,  from Lemma \ref{key}, we know that
\begin{align*}
 BT(\eta) &= (S(\eta), \eta)_{\partial M} \text{ if $ \eta \in D$ }, \\
      &=(T(\eta), \eta)_{\partial M} \text{ if $ \eta \in N$}.
\end{align*}   Here, the operators $S$ and $T$ depend only on the second fundamental form of $\partial M$.  By compactness, the second fundamental form is bounded on $\partial M$. Therefore,
\[
 |(S(\eta), \eta)|\leq \|\eta\|_{L^2(\partial M)} \|S(\eta)\|_{L^2(\partial M)} \leq C_s\|\eta\|_{L^2(\partial M)},\]
\[ 
  |(T(\eta), \eta)|\leq \|\eta\|_{L^2(\partial M)} \|S(\eta)\|_{L^2(\partial M)} \leq C_s\|\eta\|_{L^2(\partial M)}.
\] We apply Ehrling's inequality and get the estimate
\[
\|\eta\|^2_{L^2(\partial M)} \leq \epsilon \|\eta\|^2_{1}+C_{\epsilon}\|\eta\|^2_{0}.
\] Therefore
\begin{itemize}
\item $|BT(\eta)| \leq \epsilon \|\eta\|^2_{1}+C_{s,\epsilon}\|\eta\|^2_{0}$ when $\eta \in D$ or $\eta \in N$;
\item $|BT(\mathcal{J}\eta)| \leq  \epsilon \|\mathcal{J}\eta\|^2_{1}+C_{s,\epsilon}\|\mathcal{J}\eta\|^2_{0}\\
\leq  \epsilon\left( \|\eta\|^2_{1}+C_{J}\|\eta\|^2_{0}\right)+C_{s,\epsilon}\|\eta\|^2_{0}$ when $\mathcal{J}\eta \in D$ or $\mathcal{J}\eta \in N$. Here $C_{J}$ is a constant depending only on the derivative of $\mathcal{J}$.
\end{itemize}
 Thus, when $n-k>0$, we obtain
\[
C(\nabla \eta, \nabla \eta)\leq \left(D_{J}(\eta, \eta)+\|\eta\|^2_{0}\right)
\] for any $\eta\in P^k$satisfying both of the following boundary conditions:
\begin{equation*}
\left\{\begin{aligned}
\eta &\in D\, \text{or}\, \eta \in N,\\
\mathcal{J}\eta &\in D\, \text{or}\, \mathcal{J}\eta \in N.
\end{aligned}\right.
\end{equation*} These boundary conditions are equivalent to the following:
\begin{equation*}B_1:
\left\{\begin{aligned}
\eta &\in D\, \text{or}\, \eta \in N,\\
\eta &\in JD\, \text{or}\,\eta \in JN.
\end{aligned}\right.
\end{equation*} 

Since $\|\eta\|^2_{1}=(\nabla\eta, \nabla\eta)$, we get 
\[
C\|\eta\|^2_{1} \leq  \left(D_{J}(\eta, \eta)+\|\eta\|^2_{0}\right)
\] for any $\eta \in P^k$ satisfying the boundary condition $B_1$. Here, the constant $C$ depends only on $(\omega, J,g)$.
\end{proof}

\begin{proof}[Proof of Gaffney's inequality for $\Delta_{+}$ and $\Delta_{-}$]
For any $\eta \in P^{k}$ with $k<n$, there is
\begin{align*}
D_{\partial_+} (\eta, \eta)&=D_J(\eta, \eta)+ \frac{1}{n-k+1}(d^{\Lambda\ast}\eta, d^{\Lambda\ast}\eta)\\
&\geq D_J(\eta, \eta).
\end{align*} Since
\begin{equation*}
\eta \in D \Rightarrow \eta \in JN,\,  \eta\in JD\Rightarrow \eta \in N,
\end{equation*} $\eta$ satisfies the boundary condition $B_1$ in Lemma \ref{cor1} above when it satisfies the boundary condition $D$ or $JD$. Gaffney's inequality for $\Delta_+$ thus follows. The case for $D_{\partial_{-}}$ is an immediate consequence of the case of $D_{\partial_{+}}$ by conjugation.
\end{proof}

\end{document}